\documentclass[12pt, leqno]{article}
\usepackage{amssymb}
\def\q \m#1#2{{\raise1pt\hbox{$#1$}\kern-1pt\big/
               \kern-1pt\raise-1pt\hbox{$#2$}}}

\def\bJ{{\rm \bf J}}
\def\bN{{\rm \bf N}}
\def\bR{{\rm \bf R}}
\def\bZ{{\rm \bf Z}}
\def\bQ{{\rm \bf Q}}
\def\bC{{\rm {\bf C}}}
\def\ch{{\rm  ch}}

\def\sR{{ \rm \scriptsize  \bf R}}
\def\sZ{{ \rm \scriptsize  \bf Z}}

\def\sm{{ \rm \scriptsize   mod}}
\def\md{{ \rm    mod}}

\newfam\msbfam
\font\twelmsb=msbm10 at 12pt
\font\tenmsb=msbm10 at 10 pt
\font\sevenmsb=msbm10 at 7pt

\textfont\msbfam=\twelmsb
\textfont\msbfam=\tenmsb
\scriptfont\msbfam=\sevenmsb

\newtheorem{thm}{Theorem}[section]
\newtheorem{lemma}{Lemma}[section]

 \newtheorem{Rem}{Remark}[section]
 \newtheorem{prop}{Proposition}[section]

\begin{document}

\renewcommand{\theequation}{\thesection.\arabic{equation}}
\setcounter{equation}{0}

\centerline{\Large {\bf Rigidity and Vanishing Theorems in $K$-Theory II }}

\vskip 10mm  
\centerline{\bf Kefeng   LIU\footnote{Partially supported by 
the Sloan Fellowship and an NSF grant.}, 
Xiaonan MA\footnote{Partially supported by SFB 288.}
and Weiping ZHANG\footnote{Partially supported by NSFC, MOEC and the Qiu Shi Foundation.}}
\vskip 8mm

{\bf Abstract.} We extend our family rigidity and vanishing 
theorems in [{\bf LiuMaZ}] to the Spin$^c$ case. 
In particular, we prove a $K$-theory version of the main results of [{\bf H}], 
[{\bf Liu1}, Theorem B] for a family of almost complex manifolds.\\

{\bf 0 Introduction}.
Let $M,\ B$ be two compact smooth manifolds, 
and $\pi: M\to B$ be a  smooth fibration with compact fibre $X$.
Assume that a compact Lie group $G$ acts fiberwise on $M$, that is, the action
preserves each fiber of $\pi$. Let $P$ be a family of $G$-equivariant elliptic operators
along the fiber $X$. Then the family index of $P$, 
${\rm Ind} (P)$, 
is a well-defined element in $K(B)$ (cf. [{\bf AS}])
and is a virtual $G$-representation (cf. [{\bf LiuMa1}]).  
We denote by $({\rm Ind} (P))^G\in K(B)$ the $G$-invariant part of 
${\rm Ind} (P)$.

A family of elliptic operator $P$ is said to be
{\em rigid on the equivariant Chern character level }
with respect to this $G$-action, if the equivariant Chern character 
$\ch_g ({\rm Ind} (P)) \in H^* (B)$ is independent of $g \in G$. 
If $\ch_g ({\rm Ind} (P))$ is identically zero for any $g$,
 then we say $P$ has 
{\em vanishing property on the equivariant Chern character level}. 
More generally, we say that $P$ is 
{\em rigid on the equivariant $K$-theory level}, if 
${\rm Ind} (P) = ({\rm Ind} (P))^G$.  
If this index is identically zero
in $K_{G}(B)$, then we say that $P$ has 
{\em vanishing property on the equivariant $K$-theory level}. 
To study rigidity and 
vanishing, we only need  to restrict to the case where $G=S^1$. 
{}From now on we assume $G=S^1$. 

As was remarked in [{\bf LiuMaZ}], the rigidity and vanishing properties on 
 the $K$-theory level are more subtle than that on the Chern character level. 
The reason is that the 
Chern character can kill the torsion elements involved in the index bundle.

In [{\bf LiuMaZ}], we proved several rigidity and vanishing theorems 
on the equivariant $K$-theory level for elliptic genera. 
In this paper, we apply the method in [{\bf LiuMaZ}] 
to prove   rigidity and vanishing 
 theorems on the equivariant $K$-theory level for Spin$^c$ manifolds, as well 
as for almost complex manifolds.
To prove the  main results of this paper, to be stated in Section 2.1, 
we will introduce some 
shift operators on certain  vector bundles over the fixed point 
set of the circle action, 
and compare the index bundles after the shift operation. 
Then we get a recursive  relation of these index bundles which will 
in turn lead us to the final result (cf. [{\bf LiuMaZ}]). 

Let us state some of our main results in this paper
more explicitly. As was remarked in [{\bf LiuMaZ}], our
method is inspired by the ideas of Taubes [{\bf T}] and Bott-Taubes [{\bf BT}].

For a complex (resp. real) vector bundle $E$ over $M$, let
\begin{eqnarray}\begin{array}{l}
{\rm Sym}_t (E) = 1 + t E + t^2 {\rm Sym}^2 E + \cdots ,\\
\Lambda_t (E) = 1 + tE + t^2 \Lambda^2 E + \cdots
\end{array}\end{eqnarray}
be the symmetric and exterior power   operations of $E$ 
(resp. $E\otimes_\sR {\bf C}$)
in $K(M)[[t]]$ respectively.

We assume that $TX$ has an $S^1$-invariant almost complex structure $J$.
Then we can construct  canonically the Spin$^c$ Dirac operator $D^X$ 
on $\Lambda^* (T^{(0,1)*}X)$ along the fiber $X$. 
Let $W$ be an $S^1$-equivariant complex vector bundle over $M$. 
We denote by $K_W= \det W$ and $K_X= \det (T^{(1,0)} X)$ the determinant line bundles 
of $W$ and $T^{(1,0)}X$ respectively. Let 
\begin{eqnarray}\begin{array}{l}
Q_1 (W)= \otimes_{n=0}^{\infty} \Lambda_{-q^n} (\overline{W})
\otimes \bigotimes_{n=1}^{\infty} \Lambda_{-q^n} ({W}).
\end{array}\end{eqnarray}
For $N\in \bN^*$, let $y=e^{2 \pi i /N} \in \bC$. 
Let $G_y$ be the multiplicative group 
generated by $y$.  Following Witten [{\bf W}], we consider the fiberwise action $G_y$
 on $W$ and $\overline{W}$ by sending $y\in G_y$ to $y$ on $W$ and $y^{-1}$ 
on $\overline {W}$. Then $G_y$ acts naturally on $Q_1(W)$. 
We define $Q_1( T^{(1,0)} X)$ and the action $G_y$ on it in the above way.

The following theorem generalizes the result in [{\bf H}] to the family case.
\begin{thm} Assume  $c_1(T^{(1,0)}X) =0 \ \md (N) $, the family of $G_y \times S^1$ 
equivariant Spin$^c$ Dirac operators 
$D^X \otimes_{n=1}^{\infty} {\rm Sym}_{q^n} 
(TX\otimes_{\sR} \bC) \otimes Q_1 (T^{(1,0)} X)$   
is rigid on the equivariant K-theory level, for the $S^1$ action.
\end{thm}

The following family rigidity and vanishing theorem generalizes 
[{\bf Liu1}, Theorem B] to the family case.
\begin{thm} 
Assume $\omega_2(TX-W)_{S^1}=0,\ {1 \over 2} p_1(TX-W)_{S^1}= 
e{\cdot} {\bar \pi}^* u^2\ (e\in \bZ)$ in $H^*_{S^1}(M, \bZ)$, 
and $c_1(W) =0 \ \md (N) $. Consider the family of 
$G_y \times S^1$ equivariant Spin$^c$ Dirac operators 
$$ D^X \otimes (K_W\otimes K_X^{-1})^{1/2} \otimes_{n=1}^{\infty} 
{\rm Sym}_{q^n} (TX\otimes_{\sR} \bC) \otimes Q_1 (W). $$

i) If $e=0$, then these operators are rigid on the equivariant K-theory level
 for the $S^1$ action.

ii) If $e<0$, then the index bundles of these operators are zero in 
$K_{G_y \times S^1} (B)$. In particular, these index bundles are zero
 in $K_{G_y}(B)$. 
\end{thm}
We refer to  Section 2 for more details on the notation in Theorem 0.2.
Actually, our main result, Theorem 2.2, holds on a family of Spin$^c$-manifolds 
with Theorem 0.2 being one of its special cases.

This paper is organized as follows. In Section 1, we recall 
 a $K$-Theory version 
of the equivariant family index theorem for the circle action case 
[{\bf LiuMaZ}, Theorem 1.2]. 
As an immediate corollary, we get a $K$-theory version of the vanishing 
theorem of Hattori  for a family of almost complex manifolds.
In Section 2, we prove  the rigidity and vanishing theorem for elliptic genera 
in the Spin$^c$ case, on the equivariant $K$-theory level. 
The proof of the main results in Section 2 is 
base on two intermediate results which will be proved in Sections 3 and 4
respectively.

$\ $

{\bf Acknowledgements.} Part of this work was done while the authors were 
visiting the
Morningside Center for Mathematics in Beijing during the summer of 1999. 
The authors would like to 
thank the Morningside Center for hospitality. The second author would also like to 
thank the Nankai Institute of Mathematics for hospitality.

\newpage 

\section{ \normalsize  A $K$-theory version of the  equivariant family 
index theorem}
\setcounter{equation}{0}

In this section, we  recall  a $K$-theory version 
of the equivariant family  index theorem
 [{\bf LiuMaZ}, Theorem 1.2] for  $S^1$-actions, which will play
a crucial role in the following sections. 

This section is organized as follows: 
In Section 1.1, we recall the $K$-theory version of the equivariant family
index theorem for $S^1$-actions on a family of Spin$^c$ manifolds. 
In Section 1.2,  as a simple application of Theorem 1.1, we obtain a $K$-theory 
version of the vanishing theorem of Hattori [{\bf Ha}]
for the case of almost complex manifolds.

\subsection{\normalsize A $K$-theory version of the equivariant family index theorem}

Let $M,\ B$ be two compact manifolds,
let $\pi: M\to B$ be a fibration with compact fibre $X$ such that 
$\dim X = 2 l$
 and that $S^1$ acts fiberwise on $M$.  Let $h^{TX}$ be a metric on $TX$.
We assume that $TX$ is oriented.
Let $(W, h^W)$ be a Hermitian complex  vector bundle over $M$.

Let $V$ be a $2p$ dimensional oriented real vector bundle   over $M$. 
Let $L$ be a complex line bundle over $M$ with the property that the vector 
bundle $U= TX \oplus V$ obeys $\omega_2(U) = c_1(L) \ {\rm mod} \ ( 2)$. 
Then the vector 
bundle $U$ has a Spin$^c$-structure.  Let $h^V,\ h^L$ be the corresponding metrics on $V,\ L$.
 Let $S(U,L)$ be the fundamental complex spinor bundle for $(U, L)$ 
[{\bf LaM}, Appendix D.9] which  locally may be written  as 
\begin{eqnarray}
S(U,L)= S_0(U) \otimes L^{1/2},
\end{eqnarray}
where $S_0(U)$ is the fundamental spinor bundle 
for the (possibly non-existent) spin structure on $U$, 
and where $L^{1/2}$ is the 
(possibly non-existent) square root of $L$.

Assume that the $S^1$-action on $M$ lifts to  $V$, $L$ and $W$, 
and assume the metrics  $h^{TX},\ h^V,\ h^L,\ h^W$ are $S^1$-invariant.
 Also assume  that  the $S^1$-actions on $TX,\ V,\ L$ lift to  $S(U,L)$.  

Let $\nabla^{TX}$ be the Levi-Civita connection on $(TX, h^{TX})$ along 
the fibre $X$. Let $\nabla^V$, $\nabla^L$ and $ \nabla^W$ be the $S^1$-invariant
 and metric-compatible connections on $(V, h^V)$, $(L, h^L)$ and  $(W, h^W)$
respectively. 
Let $\nabla^{S(U, L)}$ be the Hermitian connection on $S(U, L)$ induced by 
$\nabla^{TX} \oplus \nabla^V$ and  $\nabla^L$ 
(cf. [{\bf LaM}, Appendix D], [{\bf LiuMaZ}, \S 1.1]).
Let $\nabla^{S(U,L)\otimes W }$ be the tensor product connection on 
$S(U,L) \otimes W$ induced by $\nabla^{S(U,L)}$ and  $\nabla^W$,
\begin{eqnarray}
\nabla^{S(U,L) \otimes W}= \nabla^{S(U,L) } \otimes 1 + 1 \otimes \nabla^W.
\end{eqnarray}

Let $\{e_i\}_{i=1}^{2l}$ (resp. $\{ f_j \}_{j=1}^{2p}$)
 be  an oriented orthonormal 
basis of $(TX, h^{TX})$ (resp. $ (V, h^V)$). We denote by $c(\cdot)$ the Clifford 
action of $TX\oplus V$ on $S(U,L)$.  Let $D^X\otimes W$ be the family 
Spin$^c$-Dirac operator on the fiber $X$ defined by
\begin{eqnarray}\label{dirac1}
D^X  \otimes W= \sum_{i=1}^{2l} c(e_i) \nabla^{S(U,L)\otimes W}_{e_i} .
\end{eqnarray}

There are two canonical  ways to consider $S(U,L)$ as a $\bZ_2$-graded vector 
bundle. Let 
\begin{eqnarray}\begin{array}{l}
\tau_s = i^{l } c(e_1) \cdots c(e_{2l}),\\
\tau_e = i^{l+p} c(e_1)\cdots c(e_{2l}) c(f_1) \cdots c(f_{2p})
\end{array}\end{eqnarray}
be two involutions of $S(U, L)$. Then $\tau_s^2= \tau_e^2=1$. We decompose 
$S(U,L)= S^+(U,L) \oplus S^-(U,L)$ corresponding to $\tau_s$ (resp. $\tau_e$) 
such that $\tau_s|_{S^{\pm}(U,L)} = \pm 1$ 
(resp. $\tau_e|_{S^{\pm}(U,L)} = \pm 1$).

For $\tau = \tau_s $ or $\tau_e$, 
by [{\bf LiuMa1}, Proposition 1.1], the index bundle 
${\rm Ind}_\tau (D^X)$ over $B$
is well-defined in the equivariant $K$-group $K_{S^1}(B)$.

Let $F=\{F_{\alpha}\} $ be the fixed point set of the circle action on $M$. 
Then $\pi: F_{\alpha}\to B$ (resp. $\pi: F\to B$) is a 
smooth fibration with fibre $Y_{\alpha}$ (resp. $Y$). 
Let $\widetilde{\pi}: N\to F$ denote the normal bundle to $F$ in $M$. 
Then $N= TX/TY$. We identify $N$ as the orthogonal complement 
of $TY$ in $TX_{|F}$.  
Let $h^{TY},\ h^N$ be the corresponding metrics on $TY$ and $ N$ induced by $h^{TX}$.
Then, we have the following $S^1$-equivariant decomposition  of $TX$ over $F$,
\begin{eqnarray}
TX_{|F} = N_{m_1} \oplus \cdots \oplus N_{m_l} \oplus TY, \nonumber
\end{eqnarray}
where each $N_{\gamma}$ is a complex vector bundle such that $g\in S^1$ acts 
on it by $g^{{\gamma} }$.  To simplify the notation, we will write simply that
\begin{eqnarray}\label{dirac6}
TX_{|F} = \oplus_{v\neq 0} N_v \oplus TY, 
\end{eqnarray}
where $N_v$ is a complex vector bundle such that $g\in S^1$ acts on it 
by $g^v$ with  $v\in \bZ^*$.
Clearly, $N= \oplus_{v\neq 0} N_v$.
We will denote by $N$ a complex vector bundle, and $N_{\sR}$ the underlying 
real vector bundle of $N$.

Similarly let
\begin{eqnarray}\label{dirac7}
W_{|F} = \oplus_v W_{v}
\end{eqnarray}
be the $S^1$-equivariant decomposition of the restriction of $W$ over $F$. 
Here $W_v \ (v\in \bZ)$ is a complex vector bundle over $F$ on which 
$g\in S^1$ acts by $g^v$. 

We also have the following $S^1$-equivariant decomposition of $V$ restricted to $F$,
\begin{eqnarray}\label{dirac2}
V_{|F}= \oplus_{v\neq 0} V_v \oplus V_{0}^{\sR},
\end{eqnarray}
where $V_v$ is a complex vector bundle such that $ g$ acts on it by $g^v$, 
and $V_0^{\sR}$ is the real subbundle of $V$ such that $S^1$ acts  
 as identity. For $v\neq 0$, let $V_{v,\bR}$ denote the underlying real 
vector bundle of $V_v$.
Denote  by  $2p'= \dim V_0^{\sR}$ and $ 2l'= \dim Y$. 

Let us write 
\begin{eqnarray} \label{dirac3}
L_F = L \otimes \Big (\bigotimes _{v\neq 0} \det N_v \bigotimes _{v\neq 0} 
\det V_v
\Big )^{-1}.
\end{eqnarray}
Then $TY \oplus V_0^{\sR}$ has a Spin$^c$ structure as
$\omega_2(TY\oplus V_0^{\sR}) = c_1(L_F) \ {\rm mod} \  (2)$. 
Let $S(TY \oplus V_0^{\sR}, L_F)$ be the fundamental spinor bundle for 
$(TY\oplus V_0^{\sR}, L_F)$ [{\bf LaM}, Appendix D, pp. 397].  

Let $D^Y,\ D^{Y_\alpha}$ be the families of Spin$^c$ Dirac operators  acting on
$ S(TY \oplus V_0^{\sR}, L_F)$ over $F,\ F_\alpha$ as (\ref{dirac1}). 
If $R$ is an Hermitian complex vector bundle equipped with an Hermitian
connection over $F$, let $D^Y\otimes R,\ D^{Y_\alpha}\otimes R$ denote 
the twisted Spin$^c$ Dirac operators on $ S(TY \oplus V_0^{\sR}, L_F) \otimes R$ and on $ S(TY_\alpha \oplus V_0^{\sR}, L_F) \otimes R$ respectively. 

Recall that $N_{v, \sR}$ and $V_{v,\sR}$  are canonically oriented by their
 complex structures. The decompositions (\ref{dirac6}), ({\ref{dirac2}) 
induce  the orientations 
on $TY$ and $V_0^{\sR}$ respectively. Let $\{e_i\}_{i=1}^{2l'}$, $\{ f_j \}_{j=1}^{2p'}$ 
be  the corresponding oriented orthonormal basis of 
$(TY, h^{TY})$ and $ (V_0^{\sR}, h^{V_0^{\sR}})$. 
There are two canonical ways to consider 
$ S(TY \oplus V_0^{\sR}, L_F) $ as a  $\bZ_2$-graded vector bundle. Let 
\begin{eqnarray}\begin{array}{l}
\tau_s = i^{l' } c(e_1) \cdots c(e_{2l'}),\\
\tau_e = i^{l'+p'} c(e_1)\cdots c(e_{2l'}) c(f_1) 
\cdots c(f_{2p'})
\end{array}\end{eqnarray}
be two involutions of $S(TY \oplus V_0^{\sR}, L_F)$. 
Then $\tau_s^2= \tau_e^2=1$. We decompose 
$S(TY \oplus V_0^{\sR}, L_F)= S^+(TY \oplus V_0^{\sR}, L_F)$ 
$\oplus S^-(TY \oplus V_0^{\sR}, L_F)$ corresponding to $\tau_s$ (resp. $\tau_e$) 
such that $\tau_s|_{S^{\pm}(TY \oplus V_0^{\sR}, L_F)} = \pm 1$ 
(resp. $\tau_e|_{S^{\pm}(TY \oplus V_0^{\sR}, L_F)} = \pm 1$).

Upon restriction to $F$, one has the following isomorphism of $\bZ_2$-graded 
Clifford modules over $F$,
\begin{eqnarray} \label{dirac4}
S(U, L) \simeq S(TY \oplus V_0^{\sR}, L_F) 
\widehat{\bigotimes _{v\neq 0}} \Lambda  N_v \widehat{\bigotimes _{v \neq 0}} \Lambda  V_v.
\end{eqnarray}
We denote by ${\rm Ind}_{\tau_s}$, ${\rm Ind} _{\tau_e}$ the index bundles 
corresponding to the involutions $\tau_s,\ \tau_e$ respectively.

Let $S^1$ act on $L$ by sending $g\in S^1$ to $g^{l_c}$ $(l_c\in \bZ)$ on $F$. 
Then $l_c$ is locally constant on $F$.
We define the following elements in $K(F) [[q^{1/2}]]$,
\begin{eqnarray}  \label{dirac7}\quad \begin{array}{l} 
R_{\pm}(q) = q^{{1 \over 2}\Sigma_v |v| \dim N_v 
- {1 \over 2} \Sigma_v v \dim V_v +{1 \over 2} l_c}
 \otimes_{0<v}\Big ( {\rm Sym}_{q^v} (N_v) \otimes \det N_v\Big )\\
\hspace*{20mm}\otimes_{v<0} {\rm Sym}_{q^{-v}} (\overline{N}_v) 
\otimes _{v\neq 0} 
\Lambda_{\pm q^v} (V_v) \otimes _v q^v W_v = \sum_n R_{\pm, n} q^n,\\
R'_{\pm}(q) = q^{-{1 \over 2}\Sigma_v |v| \dim N_v 
- {1 \over 2} \Sigma_v v \dim V_v +{1 \over 2} l_c} 
\otimes_{0<v} {\rm Sym}_{q^{-v}} (\overline{N}_v) \\
\hspace*{10mm}\otimes_{v<0}\Big ( {\rm Sym}_{q^{v}} ({N}_v) 
\otimes \det N_v \Big )\otimes _{v\neq 0} \Lambda_{\pm q^v} (V_v) 
\otimes _v q^v W_v = \sum_n R'_{\pm, n} q^n.
\end{array}\end{eqnarray}

The following result was proved in [{\bf LiuMaZ}, Theorem 1.2]:

\begin{thm} For $n\in \bZ$, we have the following identity in $K(B)$,
\begin{eqnarray}\begin{array}{l} 
{\rm Ind}_{\tau_s} (D^X \otimes W, n) = 
\sum_\alpha (-1)^{\Sigma_{0<v} \dim N_v} 
{\rm Ind} _{\tau_s} (D^{Y_\alpha} \otimes R_{+, n})  \\
\hspace*{30mm}= \sum_\alpha (-1)^{\Sigma_{v<0} \dim N_v} 
{\rm Ind} _{\tau_s} (D^{Y_\alpha} \otimes R'_{+, n}) ,\\
 {\rm Ind}_{\tau_e} (D^X\otimes W, n) 
= \sum_\alpha (-1)^{\Sigma_{0<v} \dim N_v} 
{\rm Ind} _{\tau_e} (D^{Y_\alpha} \otimes R_{-, n})  \\
\hspace*{30mm}= \sum_\alpha (-1)^{\Sigma_{v<0} \dim N_v} 
{\rm Ind} _{\tau_e} (D^{Y_\alpha} \otimes R'_{-, n}).
\end{array}\end{eqnarray}
\end{thm}

{\bf Remark 1.1.} If $TX$ has an $S^1$-equivariant Spin structure, 
by setting $V=0, L=\bC$, we get [{\bf LiuMaZ}, Theorem 1.1].

\subsection{ \normalsize $K$-theory version of the vanishing theorem of Hattori}

In this subsection, we assume that $TX$ has an $S^1$-equivariant almost 
complex structure $J$. Then one has the canonical splitting
\begin{eqnarray}
TX \otimes_{\sR} \bC = T^{(1,0)}X \oplus T^{(0,1)}X ,
\end{eqnarray}
where
\begin{eqnarray}\begin{array}{l} 
T^{(1,0)}X = \{ z\in TX \otimes_{\sR} \bC, Jz = \sqrt{-1} z\},\\
T^{(0,1)}X = \{ z\in TX \otimes_{\sR} \bC, Jz = -\sqrt{-1} z\}.
\end{array}\nonumber
\end{eqnarray}
Let $K_X = \det (T^{(1,0)}X )$ be the determinant line bundle of $T^{(1,0)}X$ 
over $M$. Then the complex spinor bundle $S(TX, K_X)$  for $(TX, K_X)$ is 
$\Lambda (T^{(0,1)*}X)$. In this case, the almost complex structure $J$ 
on $TX$ induces an almost complex structure on $TY$.
Then we can rewrite (\ref{dirac6}) as,
\begin{eqnarray}
T^{(1,0)}X = \oplus_{v \neq 0} N_v \oplus T^{(1,0)}Y,
\end{eqnarray}
where $N_v$ are complex vector subbundles of $T^{(1,0)}X $ on which 
$g\in S^1$ acts by multiplication by $g^v$.

We suppose that $c_1(T^{(1,0)}X)= 0 \ \md (N)$ $(N\in \bZ, N\geq 2)$. 
Then the complex line bundle $K_X^{1/N}$ is well defined over $M$.
After replacing the $S^1$ action by its $N$-fold action, we can always 
assume that $S^1$ acts on $K_X^{1/N}$.
For $s\in \bZ$, let $D^X \otimes K_X^{s/N}$ be the twisted Dirac operator on 
$\Lambda (T^{(0,1)*}X) \otimes K_X^{s/N}$ defined as in (\ref{dirac1}). 

The following result  generalizes the main result of
[{\bf Ha}] to the family case.

\begin{thm}  We assume that $M$ is connected and that the $S^1$ action is nontrivial.
If $c_1(T^{(1,0)}X)= 0 \ \md (N)$ $(N\in \bZ, N\geq 2)$, then for 
$s\in \bZ,\ -N < s < 0$,
\begin{eqnarray}
{\rm Ind}( D^X \otimes K_X^{s/N}) =0 \ {\rm in } \    K_{S^1}(B).
\end{eqnarray}
\end{thm}

{\em Proof} : Consider $R_+(q),\ R_+'(q)$ of (\ref{dirac7}) with $V=0, W= K_X^{s/N}$. 
We know 
\begin{eqnarray} \label{dirac8}\begin{array}{l} 
R_{+,n}=0\  {\rm if} \    n<a_1= {\rm inf}_{\alpha}({1\over 2} 
\sum_v |v| \dim N_v + ({1 \over 2}+ {s \over N}) \sum_v v \dim N_v) ,\\
R'_{+,n}=0\  {\rm if} \    n>a_2 ={\rm sup}_{\alpha}(-{1\over 2} 
\sum_v |v| \dim N_v + ({1 \over 2}+ {s \over N}) \sum_v v \dim N_v) .
\end{array}
\end{eqnarray}
As $-N < s < 0$, by (\ref{dirac8}), we know that $a_1 \geq 0,\ a_2\leq 0$,
with $a_1$ or $a_2$ equal to zero iff $\sum_v |v| \dim N_v =0$ for all $\alpha$,
 which means that the $S^1$ action does not have fixed points. 

{} From Theorem 1.1 (cf. [{\bf Z}, Theorem A.1]) and the above discussion, 
we get Theorem 1.2. \hfill $\blacksquare$

$\ $

{\bf Remark 1.2.} From  the proof of Theorem 1.2, one also deduces
that $D^X \otimes K_X^{-1},\ D^X$ 
are rigid  on the equivariant $K$-theory level (cf. [{\bf Z}, (2.17)]).
\newpage

\section{ \normalsize Rigidity and vanishing theorems in K-Theory}
\setcounter{equation}{0}

The purpose of this section is to establish the main results of this paper: 
the  rigidity and vanishing theorems on the equivariant $K$-theory level
for a family of Spin$^c$ manifolds.
The results in this section 
refine some of the results in [{\bf LiuMa2}] to the $K$-theory level.

 This section is organized as follows: 
In Section 2.1, we state our main results, the rigidity and vanishing 
theorems on  the equivariant $K$-theory level 
for a family of Spin$^c$ manifolds.  In Section 2.2, we state two 
intermediate results which will be used to prove our main results stated in 
Section 2.1. 
In Section 2.3, we prove the family rigidity and vanishing theorems.

Throughout this section, we keep the  notations of Section 1.1.

\subsection{ \normalsize Family rigidity and vanishing Theorem}

Let $\pi: M\to B$ be a  fibration of compact manifolds with fiber $X$ 
and $\dim X= 2l$. We assume that  $S^1$ acts fiberwise on $M$,
and $TX$ has an $S^1$-invariant Spin$^c$ structure. 
 Let $V$ be an even dimensional real vector bundle  over $M$. 
We  assume that $V$ has an $S^1$-invariant spin structure.
Let $W$ be an $S^1$-equivariant  complex vector bundle of rank $r$ over $M$.
 Let $K_W = \det (W)$ be the determinant line bundle of $W$.

Let $K_X$ be the $S^1$-equivariant complex line bundle over $M$ which is induced
by the $S^1$-invariant Spin$^c$ structure of $TX$.
 Its equivariant first Chern class $c_1(K_X)_{S^1}$ may also be written as 
$c_1(TX)_{S^1}$. 

Let $S(TX, K_X)$ be the complex spinor bundle of $(TX, K_X)$ as in Section 1.1.
Let $S(V) = S^+ (V) \oplus S ^- (V) $ be the spinor bundle of $V$.

We  define the following elements in $K(M) [[q^{1/2}]]$:
\begin{eqnarray}\begin{array}{l}
Q_1(W)= \bigotimes_{n=0}^{\infty} \Lambda_{-q^n} (\overline{W})
\otimes \bigotimes_{n=1}^{\infty} \Lambda_{-q^n} ({W}),\\
R_1(V)=(S^+ (V)+S^- (V)) \otimes_{n=1}^\infty \Lambda_{q^n} (V) ,\\
R_2(V)=(S^+ (V)- S^- (V))\otimes_{n=1}^\infty \Lambda_{-q^n} (V), \\
R_3(V)=\otimes_{n=1}^\infty \Lambda_{-q^{n-1/2}} (V), \\
R_4(V)=\otimes_{n=1}^\infty \Lambda_{q^{n-1/2}} (V).
\end{array}
\end{eqnarray}
 For $N\in \bN^*$, let $y=e^{2 \pi i /N}\in \bC$. 
Let $G_y$ be the multiplicative group 
generated by $y$.  Following Witten [{\bf W}], we consider the fiberwise action $G_y$
 on $W$ and $\overline{W}$ by sending $y\in G_y$ to $y$ on $W$ and $y^{-1}$ 
on $\overline {W}$. Then $G_y$ acts naturally on $Q_1(W)$. 

Recall that the equivariant cohomology group $H^*_{S^1} (M, \bZ)$ 
of $M$ is defined by
\begin{eqnarray}
H^*_{S^1} (M, \bZ)= H^*(M \times_{S^1} ES^1, \bZ),
\end{eqnarray}
where $ES^1$ is the usual universal $S^1$-principal bundle over the 
classifying space $BS^1$  of $S^1$.
So $H^*_{S^1} (M, \bZ)$ is a module over $H^*(BS^1, \bZ)$ induced by the 
projection $\overline{\pi} : M\times _{S^1} ES^1\to BS^1$. 
Let $p_1(V)_{S^1},\ p_1(TX)_{S^1} \in H^*_{S^1} (M, \bZ)$ be the $S^1$-equivariant
 first Pontrjagin classes of $V$ and $TX$ respectively. 
As $V\times_{S^1}ES^1$ is spin over $M\times_{S^1}ES^1$, 
one knows that ${1\over 2}p_1(V)_{S^1}$ is well-defined in $H^*_{S^1}(M,\bZ)$
(cf. [{\bf T}, pp. 456-457]).
Also recall that 
\begin{eqnarray}\label{hyp0}
H^*(BS^1, \bZ)= \bZ [[u]]
\end{eqnarray}
with $u$ a generator of degree $2$.

In the following, we denote by $D^X\otimes R$ the family of Dirac operators
acting fiberwise  
on $S(TX, K_X) \otimes R$ as was defined in Section 1.1. 

We can now state the main results of this paper as follows.

\begin{thm}\label{Th21} If $\omega_2(W)_{S^1}= \omega_2(TX)_{S^1}$,
 ${1\over 2} p_1(V+W-TX)_{S^1} = e \cdot \overline{\pi}^* u^2$
 $(n\in \bZ )$ in $H^*_{S^1} (M, \bZ)$, and 
$c_1(W)= 0 \ {\rm mod}( N)$. For  $i=1,2,3,4$, 
consider the family of $G_y \times S^1$-equivariant elliptic operators
$$ D^X \otimes (K_W\otimes K_X^{-1})^{1/2} \otimes_{n=1}^{\infty} 
{\rm Sym}_{q^n} (TX) \otimes Q_1 (W) \otimes R_i(V).$$

i) If $e=0$, then these operators are rigid on the equivariant K-theory level
 for the $S^1$ action.

ii) If $e<0$, then the index bundles of these operators are zero in 
$K_{G_y \times S^1} (B)$. In particular, these index bundles are zero
 in $K_{G_y}(B)$. 
\end{thm}
\begin{Rem}{\rm As $\omega_2(W)_{S^1}= \omega_2(TX)_{S^1}$, 
${1 \over 2}p_1(W-TX)_{S^1} \in H^*_{S^1} (M, \bZ)$ is well defined.
The condition $\omega_2(W)_{S^1}= \omega_2(TX)_{S^1}$ 
also means $c_1(K_W \otimes K_X^{-1})_{S^1}=0  \ {\rm mod} (2)$, 
by [{\bf HaY}, Corollary 1.2], the $S^1$-action on $M$ can 
be lifted to $(K_W\otimes K_X^{-1})^{1/2}$ and  is compatible 
with the $S^1$ action on $K_W\otimes K_X^{-1}$.}
\end{Rem}
\begin{Rem} {\rm If we assume $c_1(W)_{S^1}= c_1(TX)_{S^1}$ in $H^*_{S^1} (M, \bZ)$
instead of $\omega_2(W)_{S^1}= \omega_2(TX)_{S^1}$ in Theorem 2.1, then 
$K_W\otimes K_X^{-1}$ is a trivial line bundle over $M$, 
and $S^1$ acts trivially on it. In this case, Theorem 2.1 gives the 
family version of the results of 
[{\bf De}].}
\end{Rem}
\begin{Rem} {\rm The interested reader can apply our method to get various 
rigidity and vanishing theorems, for example, to get a generalization 
of Theorem1.2 for the elements [{\bf W}, (65)].}
\end{Rem}

Actually, as in [{\bf LiuMaZ}],
our proof of these theorems works under the following slightly 
weaker hypothesis. Let us first explain   some notations.

For each $n>1$, consider $\bZ_n \subset S^1$, the cyclic subgroup of order $n$.
We have the $\bZ_n$ equivariant cohomology of $M$ defined by 
$H^*_{\bZ_n}(M, \bZ) = H^*(M\times_{\bZ_n}  ES^1, \bZ)$, and there is 
a natural ``forgetful'' map $\alpha(S^1,\bZ_n): M\times_{\bZ_n} ES^1 
\to M \times_{S^1} ES^1$ which induces a pullback 
$\alpha(S^1, \bZ_n)^*: H^*_{S^1} (M, \bZ) \to H^*_{\bZ_n}(M, \bZ)$. 
The arrow which forgets the $S^1$ action altogether we denote by 
$\alpha(S^1, 1)$. Thus $\alpha(S^1, 1)^*: H^*_{S^1}(M, \bZ) \to H^*(M, \bZ)$
is induced by the inclusion of $M$ into $M\times_{S^1} ES^1$ as a fiber over 
$BS^1$. 

Finally, note that if $\bZ_n$ acts trivially on a space $Y$, then there is 
a new arrow $t^*: H^*(Y, \bZ) \to H^*_{\bZ_n} (Y, \bZ)$ induced by the 
projection $Y \times _{\bZ_n} ES^1 =Y \times B\bZ_n \stackrel{t}{\to} Y$.

We let $\bZ_{\infty}= S^1$. For each $1 < n \leq +\infty$, let $i: M(n) \to M$
 be the inclusion of the fixed point set of $\bZ_n \subset S^1$ in $M$ 
and so $i$ induces $i_{S^1}: M(n) \times_{S^1} ES^1 \to M \times _{S^1} ES^1$.  

In  the rest of this paper, we suppose  that there exists 
some integer $e\in \bZ$ such that for $1 < n \leq +\infty$,
\begin{eqnarray}\label{hyp1}
\alpha(S^1, \bZ_n)^* \circ i_{S^1} ^*  \Big ({1\over 2} p_1(V+W - TX)_{S^1} 
- e \cdot  \overline{\pi}^* u^2\Big ) \\
\hspace*{15mm}=
 t^*\circ \alpha(S^1, 1)^*\circ i_{S^1} ^*\Big ({1\over 2} p_1(V+W-TX)_{S^1}\Big ). 
\nonumber
\end{eqnarray}

\begin{Rem} {\rm The relation (\ref{hyp1}) clearly follows from the 
hypothesises  of Theorem \ref{Th21}  by pulling back and forgetting. 
Thus it is  weaker.} 
\end{Rem}

We can now state a slightly more general version of Theorem \ref{Th21}.

\begin{thm}\label{Th22} Under the hypothesis (\ref{hyp1}), we have 

i) If $e=0$, then the index bundles of the elliptic operators in 
Theorem 2.1 are rigid on the equivariant K-theory level 
for the $S^1$-action. 

ii) If $e<0$, then the index bundles of the elliptic operators in 
Theorem 2.1  are zero as  elements in $K_{G_y \times S^1}(B)$. 
In particular, these index bundles are zero in $K_{G_y}(B)$.
\end{thm}

The rest of this section is devoted to a proof of Theorem \ref{Th22}.

\subsection{\normalsize  Two intermediate results}

Let $F=\{F_ \alpha\}$ be the fixed point set of the circle action.
Then $\pi: F\to B$ is a fibration with compact fibre denoted by $Y = \{Y_ \alpha\}$.

As in [{\bf LiuMaZ}, \S 2], we may and we will assume that
\begin{eqnarray}\label{hyp2}\begin{array}{l}
TX_{|F} = TY \oplus \bigoplus_{0<v} N_{v},\\
TX \otimes_\sR \bC = TY \otimes _\sR \bC 
\bigoplus_{0<v} (N_v \oplus \overline{N}_v),
\end{array}\end{eqnarray}
where $N_v$ is the complex vector bundle on which $S^1$ acts by sending 
$g$ to $g^v$ (Here $N_v$ can be zero). We also assume that
\begin{eqnarray}\label{hyp3}\begin{array}{l}
V_{|F} = V_0^{\sR} \oplus \bigoplus_{0<v} V_v,\\
W_{|F}= \oplus_v W_v,
\end{array}\end{eqnarray}
where $V_v$, $W_v$  are complex  vector bundles on which $S^1$ acts by sending $g$ to 
$g^v$, and $V_0^{\sR}$ is a real vector bundle on which $S^1$ acts as identity.

By (\ref{hyp2}), as in (\ref{dirac4}), there is a natural 
isomorphism between the $\bZ_2$-graded $C(TX)$-Clifford modules over $F$,
\begin{eqnarray} \label{hyp20}
S(TY, K_X\otimes_{0<v} (\det N_v)^{-1}) \widehat{\otimes} _{0<v} \Lambda N_v 
\simeq S(TX, K_X)_{|F}.
\end{eqnarray}
For $R$ a complex vector bundle over $F$, let $D^Y \otimes R$, 
$D^{Y_\alpha} \otimes R$ be the twisted Spin$^c$ Dirac operator on 
$S(TY, K_X\otimes_{0<v} (\det N_v)^{-1}) \otimes R$ on $F, F_\alpha$ respectively.

On $F$, we write 
\begin{eqnarray}\label{hyp4}\begin{array}{l}
e(N) = \sum_{0<v} v^2 \dim N_v, 
\qquad d' (N) =  \sum_{0<v} v \dim N_v,\\
e(V) = \sum_{0<v} v^2 \dim V_v, \qquad d' (V) =  \sum_{0<v} v \dim V_v,\\
e(W) = \sum_{v} v^2 \dim W_v, \qquad d' (W) =  \sum_{v} v \dim W_v. 
\end{array}\end{eqnarray}
Then $e(N),\ e(V),\ e(W), \ d'(N), \ d'(V)$ and $ d'(W)$ 
are locally constant functions on $F$. 

By [{\bf H}, \S 8], we have the following property,

\begin{lemma} If $c_1(W)=0 \, \md (N)$, then $d'(W) \  {\rm mod} (N)$  
is constant on each connected component  of $M$. 
\end{lemma}

{\em Proof} : As $c_1(W)=0 \  {\rm mod} (N)$, $(K_W)^{1/N}$ is well defined. 
Consider the $N$-fold covering $S^1\to S^1$, with $\mu \to \lambda= \mu^N$, 
then $\mu$ acts on $M$ and $K_W$ through $\lambda$. 
This action can be lift to $(K_W)^{1/N}$. On $F$, $\mu$ acts on $(K_W)^{1/N}$ 
by multipication by $\mu^{d'(W)}$. However, if $\mu= \zeta= e ^{2 \pi i /N}$, 
then it operates trivially on $M$. So the action of $\zeta$  
in each fibre of $L$ is by multiplication by $\zeta^a$, 
and $a \  {\rm mod} (N)$ is constant on each connected component of $M$.

The proof of Lemma 2.1 is complete.
\hfill $\blacksquare$\\

 Let us write 
\begin{eqnarray}\label{hyp5}\begin{array}{l}
L(N) = \otimes _{0<v} (\det N_v)^v,
\qquad  L(V) = \otimes _{0<v} (\det V_v)^v,\\
L(W) = \otimes _{v\neq 0} (\det W_v)^v,\\ 
L= L(N)^{-1} \otimes L(V)\otimes L(W).
\end{array}\end{eqnarray}

We denote the Chern roots of $N_v$ by $\{ x_v ^j\}$ 
(resp. $V_v$ by $ u^j_v$ and $W_v$ by $ w^j_v$),  and the Chern roots of 
$TY \otimes_{\sR} \bC$ by $\{ \pm y_j\}$ (resp. 
$V_0=V_0^{\sR} \otimes_{\sR} \bC$ by $\{ \pm  u_0^j\}$).
 Then if we take $\bZ_\infty = S^1$ in (\ref{hyp1}), we get
\begin{eqnarray}\label{hyp6}\begin{array}{l}
{1\over 2}(\Sigma_{v,j} (u_v^j + v u)^2 +\Sigma_{v,j}(w_v^j+vu)^2
-  \Sigma_j (y_j)^2 - \Sigma_{v,j} (x_{v}^j + v u)^2 )- e u^2 \\
= {1\over 2}(\Sigma_{v,j} (u_v^j)^2+ \Sigma_{v,j}(w_v^j)^2 - \Sigma_j (y_j)^2 - \Sigma_{v,j} (x_{v}^j)^2) .
\end{array}\end{eqnarray}
By (\ref{hyp0}), (\ref{hyp6}), we get 
\begin{eqnarray}\label{hyp7}\begin{array}{l}
c_1(L) = \Sigma_{v,j} v u_{v}^j +\Sigma_{v,j} v w_{v}^j -  \Sigma_{v,j}  v x_{v}^j =0,\\
e(V)+e(W)- e(N) \\
\hspace*{5mm}= \sum_{0<v} v^2 \dim V_v+ \sum_{v} v^2 \dim W_v 
-\sum_{0<v} v^2 \dim N_v =2e,
\end{array}\end{eqnarray}
which does not depends on the connected components of $F$.
This means $L$ is a trivial complex line bundle over each component 
$F_\alpha$ of $F$, and $S^1$ acts on $L$ by sending $g$ to $g^{2e}$, 
and $G_y$ acts on $L$ by sending $y$ to $y^{d'(W)}$. 
By Lemma 2.1,  we can extend $L$ to a trivial complex line bundle over $M$, 
and we extend the $S^1$-action on it
 by sending $g$ on the canonical section $1$ of $L$ to $g^{2e} \cdot 1$, 
and $G_y$ acts on $L$ by sending $y$ to $y^{d'(W)}$.  

The line bundles in (\ref{hyp5}) will play important roles in the next two sections which
consist of the proof of Theorems \ref{Th23}, \ref{Th24} to be stated below.

In what follows, if $R(q)=\sum_{m \in {1 \over 2} \bZ} R_m q^m \in K_{S^1}(M) 
[[q^{1/2}]]$, 
we will also denote  ${\rm Ind} (D^X \otimes R_m, h)$ by
${\rm Ind} (D^X \otimes R(q), m, h)$. For $ k=1,2,3,4$, set
\begin{eqnarray}\label{hyp10}
R_{1k}= (K_W \otimes K_X^{-1})^{1 /2} \otimes Q_1(W) \otimes R_k(V).
\end{eqnarray}

We first state a result which expresses the global equivariant family index via
the family indices on the fixed point set.

\begin{prop} For $m\in {1\over 2} \bZ$, $h\in \bZ$,   $1\leq k \leq 4$,
 we have the following identity in $K_{G_y} (B)$,
\begin{eqnarray}
\begin{array}{l}
{\rm Ind} (D^X \otimes _{n=1}^\infty {\rm Sym}_{q^n} (TX) \otimes R_{1k}, m, h)\\
= \sum_\alpha (-1)^{\Sigma _{0<v} \dim N_v}{\rm Ind} 
(D^{Y_\alpha} \otimes _{n=1}^\infty {\rm Sym}_{q^n} (TX) \otimes R_{1k}\\
\hspace*{45mm}
 \otimes {\rm Sym} (\oplus_{0<v} N_v) \otimes _{0<v} \det N_v, m, h)
\end{array}
\end{eqnarray}
\end{prop}

{\em Proof} :  This follows  directly from Theorem 1.1 and 
(\ref{hyp20}).\hfill $\blacksquare$\\

For $p\in \bN$, we define the following elements in $K_{S^1}(F) [[q]]$:
\begin{eqnarray} \label{hyp8}\begin{array}{l}
{\cal F}_p (X) = \bigotimes_{0<v}\Big (  \otimes_{n=1}^{\infty}
 {\rm Sym}_{q^n} (N_v) \otimes_{n> pv} {\rm Sym}_{q^n} (\overline{N}_v)\Big )
\otimes_{n=1}^{\infty} {\rm Sym}_{q^n} (TY),\\
{\cal F}'_p (X) = \bigotimes_{\stackrel{0<v}{0\leq n \leq pv}}
\Big (  {\rm Sym}_{q^{-n}} (N_v)  \otimes \det N_v \Big ),\\
\\
{\cal F}^{-p} (X) = {\cal F}_p (X) \otimes {\cal F}'_p (X).
\end{array}\end{eqnarray}
Then, from (\ref{hyp2}), over $F$, we have 
\begin{eqnarray}
{\cal F}^0 (X) = \otimes _{n=1}^\infty {\rm Sym}_{q^n} (TX)\otimes {\rm Sym} (\oplus_{0<v} N_v) \otimes _{0<v} \det N_v.
\end{eqnarray}

We now state two intermediate results  on the relations between 
the family indices on the fixed point set. 
They will be used in the next subsection to prove Theorem \ref{Th22}.

\begin{thm}\label{Th23} For  $1\leq k \leq 4$, $h,\ p \in \bZ$, $p>0$,
$m\in {1\over 2} \bZ$,  we have the following identity in $K_{G_y}(B)$,
\begin{eqnarray}\begin{array}{l}
\sum_\alpha (-1)^{\Sigma_{0<v} \dim N_v} {\rm Ind} (D^{Y_\alpha} 
\otimes {\cal F}^0 (X)  \otimes R_{1k}, m , h)\\
=\sum_\alpha (-1)^{pd'(N)+ \Sigma_{0<v} \dim N_v} {\rm Ind} (D^{Y_\alpha} 
\otimes {\cal F}^{-p} (X)  \otimes R_{1k},\\
\hspace*{40mm} m+ {1\over 2} p^2 e(N) +{1\over 2} p d'(N), h).
 \end{array} \end{eqnarray}
\end{thm}

\begin{thm}\label{Th24} For each $\alpha$,   $1\leq k \leq 4$, $h,\ p \in \bZ$, $p>0$, 
$m\in {1\over 2} \bZ$, 
 we have the following identity in $K_{G_y}(B)$,
\begin{eqnarray}\begin{array}{l}
{\rm Ind} (D^{Y_\alpha} \otimes {\cal F}^{-p} (X)  \otimes R_{1k},
m + {1\over 2} p^2 e(N) +{1\over 2} p d'(N), h)\\
=(-1)^{p  d'(W)} {\rm Ind} (D^{Y_\alpha} \otimes {\cal F}^0 (X) 
 \otimes R_{1k} \otimes L^{-p},
m + ph + p^2 e, h ).
\end{array} \end{eqnarray}
\end{thm}

Theorem \ref{Th23} is a direct consequence of Theorem \ref{Th25} to be stated below,
which will be proved in Section 4, while Theorem \ref{Th24} will be proved in
Section 3.

To state Theorem \ref{Th25},
let $J=\{v \in \bN|$  There exists $\alpha$ such that $N_v\neq 0$ on  $F_\alpha\}$ and 
\begin{eqnarray}\label{e1}
\Phi = \{ \beta \in ]0, 1]| {\rm There \  
exists } \ v\in J \ {\rm such\ that} \ \beta v\in \bZ \}.
\end{eqnarray}
We order the elements in  $\Phi $ so that
$\Phi=\{ \beta_i| 1\leq i \leq J_0, J_0\in \bN 
\ {\rm and} \ \beta_i < \beta_{i+1}\}$.
Then for any integer $1\leq i\leq J_0$, there exist 
$p_i,\ n_i \in \bN,\ 0< p_i\leq n_i$, with $(p_i, n_i)=1$ such that 
\begin{eqnarray}\label{e2}
\beta_i = {p_i/n_i}. 
\end{eqnarray} 
Clearly, $\beta_{J_0}=1$. We also set  $p_0=0$ and  $\beta_0=0$.

For $ 1\leq j \leq J_0$, $p \in \bN^*$, we write 
\begin{eqnarray}\label{e3}\begin{array}{l}
I^p_0 = \phi, \mbox{the empty set},\\
\displaystyle{I^p_j= \{ (v,n)| v\in J, (p-1)v<n \leq  pv,   
{n \over v}= p-1 + {p_j \over n_j} \} ,}\\
\displaystyle{\overline{I}^p_j= \{ (v,n)| v\in J, (p-1)v<n \leq pv, 
{n \over v} > p-1 +  {p_j \over n_j}\}.}
\end{array}\end{eqnarray}
For $0\leq j \leq J_0$, set 
\begin{eqnarray}\label{e4}
\\
{\cal F}_{p,j} (X) = {\cal F}_p(X)\otimes {\cal F}'_{p-1} (X) 
\bigotimes_{(v,n)\in \cup_{i=1}^{j} I_i^p} \Big({\rm Sym}_{q^{-n}} ( N_v) 
\otimes  \det N_v \Big )\bigotimes_{(v,n)\in \overline{I}^p_{j}} 
{\rm Sym}_{q^{n}} ( \overline{N}_v) .\nonumber
\end{eqnarray}
Then 
\begin{eqnarray}\label{e5}\begin{array}{l}
{\cal F}_{p,0}(X) = {\cal F}^{-p+1}(X), \\
{\cal F}_{p,J_0}(X) = {\cal F}^{-p}(X).
\end{array}\end{eqnarray}

For $s\in \bR$, let $[s]$ denote the greatest integer which is less than or 
equal to the given number $s$. For $0\leq j\leq J_0$, denote by
\begin{eqnarray}\label{e6} \begin{array}{l}
e(p, \beta_j, N) = {1\over 2} \sum_{0<v} (\dim N_v )
\Big ((p-1) v + [{p_j v \over n_j}]\Big ) 
\Big ((p-1) v + [{p_j v \over n_j}]+1\Big ),\\
d'(p,\beta_j, N) =  \sum_{0<v} (\dim N_v )([{p_j v \over n_j}]+ (p-1)v).
\end{array} \end{eqnarray}
Then $e(p, \beta_j, N)$ and  $d'(p,\beta_j, N)$ 
 are locally constant functions on $F$. And 
\begin{eqnarray}\label{e7}\begin{array}{l}
e(p, \beta_0, N)={1 \over 2} (p-1)^2  e(N) + {1 \over 2} (p-1) d'(N),\\
 e(p, \beta_{J_0}, N)={1 \over 2} p^2  e(N) + {1 \over 2} p d'(N),\\
d'(p,\beta_{J_0}, N)=d'(p+1,\beta_0, N)= pd'(N).
\end{array} \end{eqnarray}

\begin{thm}\label{Th25} For   $1\leq k \leq 4$, $1\leq j \leq J_0$,
 $p \in \bN^*$, $ h\in \bZ$, $m \in {1 \over 2} \bZ$, we have the following 
identity in $K_{G_y}(B)$,
\begin{eqnarray}\label{e8}\begin{array}{l}
\sum_\alpha (-1)^{ d'(p,\beta_{j-1},N) + \Sigma_{0<v} \dim N_v} 
{\rm Ind} (D^{Y_\alpha} \otimes {\cal F}_{p,j-1} (X) 
 \otimes R_{1k},\\
\hspace*{50mm}m +e(p, \beta_{j-1}, N), h)\\
=\sum_\alpha (-1)^{d'(p,\beta_{j},N) + \Sigma_{0<v} \dim N_v} 
{\rm Ind} (D^{Y_\alpha} \otimes {\cal F}_{p, j} (X)
  \otimes R_{1k},\\
\hspace*{50mm}m + e(p, \beta_j, N) , h).
\end{array}\end{eqnarray}
\end{thm}

{\em Proof}: The proof is delayed to Section 4.\hfill $\blacksquare$\\

{\em Proof of Theorem \ref{Th23}} : From  (\ref{e5}), (\ref{e7}), 
and Theorem \ref{Th25}, 
for $1\leq k \leq 4$, $h\in \bZ$, $p\in \bN^*$ and  
$m\in  {1 \over 2} \bZ$, we have the following identity in $K_{G_y}(B)$:
\begin{eqnarray}\label{e15}
\begin{array}{l}
\sum_\alpha (-1)^{ d'(p,\beta_{J_0},N) +\Sigma_{0<v}  \dim N_v} {\rm Ind} 
(D^{Y_\alpha} \otimes {\cal F}^{-p} (X)  \otimes R_{1k},\\
\hspace*{40mm} m +{1\over 2} p^2 e(N) +{1\over 2} p d'(N), h)\\
=\sum_\alpha (-1)^{d'(p,\beta_{0},N) +\Sigma_{0<v}  \dim N_v}
{\rm Ind} (D^{Y_\alpha} \otimes {\cal F}^{-p+1} (X) \otimes 
R_{1k}, \\
\hspace*{40mm} m + {1\over 2} (p-1)^2 e(N) +{1\over 2} (p-1) d'(N), h).
\end{array}
\end{eqnarray}
{}From (\ref{e7}), (\ref{e15}),  we get Theorem  \ref{Th23}. \hfill $\blacksquare$

\subsection{ \normalsize Proof of Theorem \ref{Th22}}

As ${1\over 2} p_1(TX-W)_{S^1} \in H_{S^1}^*(M, \bZ)$ is well defined, by (\ref{hyp4}), and 
  (\ref{hyp6}), 
\begin{eqnarray}\label{e13}
d'(N)+ d'(W)=0 \ \md (2).
\end{eqnarray}

{}From Proposition 2.1, Theorems \ref{Th23}, \ref{Th24},  (\ref{e6}), (\ref{e13}),
for  
$1\leq k \leq 4$, 
$h, p\in \bZ$, $p> 0$, $m\in {1 \over 2} \bZ$,  we get the following identity 
in $K_{G_y}(B)$,
\begin{eqnarray}\label{e9}\begin{array}{l}
{\rm Ind} (D^X \otimes _{n=1}^\infty {\rm Sym}_{q^n} (TX) \otimes R_{1k},
 m, h)\\
\hspace*{10mm}
={\rm Ind} (D^X \otimes _{n=1}^\infty {\rm Sym}_{q^n} (TX) \otimes R_{1k}
\otimes L^{-p}, m', h),
\end{array}\end{eqnarray}
with 
\begin{eqnarray}\label{e10}
m'=m+ ph+ p^2 e.
\end{eqnarray}

Note that from (2.1), (\ref{hyp10}),  if $m<0$, or $m'<0$, then two side of
(\ref{e9}) are zero in $K_{G_y}(B)$. Also recall that  $y\in G_y$ acts 
on the trivial line bundle $L$ 
by sending $y$ to $y^{d'(W)}$.

i) Assume that $e=0$. Let $h\in \bZ,\ m_0\in {1 \over 2} \bZ$, $h\neq 0$ be 
fixed.
If $h>0$, 
we take $m'=m_0$, then for $p$ big enough, we get $m<0$ in (\ref{e10}).
 If $h<0$,  we take $m=m_0$,  then  for $p$ big enough, 
we get $m'<0$ in (\ref{e10}). 

 So  for $h\neq 0$, 
$m_0\in {1 \over 2} \bZ$,  $1\leq k \leq 4$, 
we get 
\begin{eqnarray}\label{e12}
{\rm Ind} (D^X \otimes _{n=1}^\infty {\rm Sym}_{q^n} (TX) 
\otimes R_{1k}, m_0, h) 
=0 \quad {\rm in } \quad  K_{G_y}(B).
\end{eqnarray}

ii) Assume that  $e <0$.  For $h\in \bZ$, $m_0\in {1 \over 2} \bZ$, 
 we take $m=m_0$, then  for $p$ big enough, 
we get $m'<0$ in (\ref{e10}), which again gives us  (\ref{e12}).

The proof of Theorem \ref{Th22} is complete. \hfill $\blacksquare$\\

{\bf Remark 2.5}:  Under the condition of Theorem 2.2 i), 
if $d'(W) \neq 0 \ {\rm mod} (N)$, 
we can't deduce these index bundles are zero in $K_{G_y}(B)$. If in addition,
$M$ is connected, by (\ref{e9}), for $1\leq k \leq 4$,
 in $K_{G_y}(B)$,  we get
\begin{eqnarray}\label{e16}\begin{array}{l}
{\rm Ind} (D^X \otimes _{n=1}^\infty {\rm Sym}_{q^n} (TX) \otimes R_{1k})\\
\hspace*{10mm}= 
{\rm Ind} (D^X \otimes _{n=1}^\infty  {\rm Sym}_{q^n} (TX) \otimes R_{1k}) 
\otimes [d'(W)].
\end{array}\end{eqnarray}
Here we denote by $[d'(W)]$ the one dimensional complex vector space  on 
which $y\in G_y$ acts by multiplication by $y^{d'(W)}$. In particular,
if $B$ is a point, by (\ref{e16}), we get the vanishing theorem analogue to 
the result of [{\bf H}, \S 10].\\

{\bf Remark 2.6}: If we replace $c_1(W)=0 \ {\rm mod} (N), y= e^{2 \pi i /N}$ 
by $c_1(W)=0, y= e ^{2 \pi  c i}$, with $c\in \bR \setminus \bQ$ 
in Theorem 2.2, then by Lemma 2.1, $d'(W)$ is constant on each connected
 component of $M$. In this case, we still have Theorem 2.2. 
In fact, we only use $c_1(W)=0 \ {\rm mod} (N)$ 
to insure the action $G_y$ on $L$ is well defined. 
So we also generalize the main result of [{\bf K}] to family case. 

\newpage

\section{ \normalsize Proof of Theorem \ref{Th24}}
\setcounter{equation}{0}

This section is organized as follows: 
In Section 3.1, we introduce some notations. In Section 3.2, we prove 
Theorem \ref{Th24} by  introducing some shift operators as in 
 [{\bf LiuMaZ}, \S 3].

Throughout this section, we keep the  notations of Section 2.

\subsection{ \normalsize Reformulation  of Theorem \ref{Th24}}

To simplify the notations, we  introduce some new notations in this 
subsection. For $n_0\in \bN^*$, we define a number operator $P$ on 
$K_{S^1}(M)[[q^{1 \over n_0}]]$ in the following way: if 
$R(q)= \oplus_{n\in {1 \over n_0} \sZ} q^n R_n \in K_{S^1}(M)[[q^{1 \over n_0}]]$,
 then $P$ acts on $R(q)$ by multiplication by $n$ on $R_n$. 
{}From now on, we simply denote ${\rm Sym} _{q^n}(TX),\ \Lambda_{q^n}(V)$ 
by ${\rm Sym} (TX_n),\ \Lambda (V_n)$ respectively. In this way, 
$P$ acts on $TX_n$, $V_n$ by multiplication by $n$, and the action $P$ on 
${\rm Sym} (TX_n),\ \Lambda (V_n)$  is naturally induced by the corresponding 
action of $P$ on $TX_n$, $V_n$. 
So the eigenspace of $P=n$ is just given by the coefficient of $q^n$ 
of the corresponding element $R(q)$. 
For $R(q)= \oplus_{n\in {1 \over n_0} \sZ} q^n R_n 
\in K_{S^1}(M)[[q^{1 \over n_0}]]$, we will also denote 
\begin{eqnarray}
{\rm Ind} (D^X \otimes R(q), m , h) = {\rm Ind} (D^X \otimes R_m, h).
\end{eqnarray}

Let $H$ be the canonical basis of ${\rm Lie} (S^1) = \bR$, 
i.e., $ \exp (t H) = {\exp}  (2 \pi i t)$ for $t\in \bR$.  
If $E$ is an $S^1$-equivariant vector bundle over $M$,
on the fixed point set $F$, let $J_H$ be the representation of 
${\rm Lie} (S^1)$ on $E|_F$. Then the weight of $S^1$ action 
on $\Gamma(F,E|_F)$ is given by the action 
\begin{eqnarray}
\bJ_H= {-1 \over 2 \pi} \sqrt{-1} J_H.
\end{eqnarray}

Recall that the $\bZ_2$ grading on 
$S(TX, K_X) \otimes _{n=1}^\infty {\rm Sym} (TX_n)$ 
(resp. $S(TY, K_X \otimes \otimes_{0<v} (\det N_v)^{-1}) \otimes {\cal F}^{-p} (X)$) 
is induced by the $\bZ_2$-grading on $S(TX, K_X)$
(resp. $S(TY, K_X \otimes \otimes_{0<v} (\det N_v)^{-1}) $).
Let
\begin{eqnarray}\begin{array}{l}
F^1_V = S(V) \bigotimes _{n=1}^\infty \Lambda (V_n) ,\\
F^2_V= \otimes_{n\in \bN + {1 \over 2}} \Lambda (V_n) ,\\
Q(W)=\otimes_{n=0}^\infty \Lambda(\overline{W}_n) \otimes_{n=1}^\infty 
\Lambda({W}_n)
\end{array}\end{eqnarray}
There are two natural $\bZ_2$ gradings on $F^1_V,\ F^2_V$ (resp. $Q(W)$). 
The first grading is induced by the
$\bZ_2$-grading of $S(V)$ and the forms of homogeneous degree in 
$\bigotimes _{n=1}^\infty \Lambda (V_n)$, 
$\otimes_{n\in \bN + {1 \over 2}} \Lambda (V_n) $ (resp. $Q(W)$). 
We define $\tau_{e|F_V^{i\pm}}= \pm 1$ (resp. $\tau_{1|Q(W)^{\pm}}= \pm 1$)
to be the involution defined by this $\bZ_2$-grading. 
The second grading is the one
for which $F^i_V$ $ (i=1,\ 2)$ are purely even, 
i.e., $F^{i+}_V=F^i_V$. We denote by $\tau_s= {\rm Id}$
  the involution 
defined by  this $\bZ_2$ grading. Then the coefficient of 
$q^n$ $ (n\in {1 \over 2} \bZ)$ in (2.1) of $R_1(V)$ or $ R_2(V)$ 
(resp. $R_3(V),\ R_4(V)$, or $ Q_1(W)$) is  exactly the  $\bZ_2$-graded vector subbundle of 
$(F^1_V, \tau_s)$ or $(F^1_V, \tau_e)$ 
(resp. $ (F^2_V, \tau_e)$, $(F^2_V, \tau_s)$ or $(Q(W), \tau_1)$), 
on which $P$ acts  by multiplication by $n$.

We  denote by $\tau_e$ (resp. by $\tau_s$) the $\bZ_2$-grading on 
$S(TX, K_X) \otimes _{n=1}^\infty {\rm Sym} (TX_n) \otimes F^k_V$ ($k=1,\ 2$)
induced by the above $\bZ_2$-gradings.
We will denote by $\tau_{e1}$ (resp. by $\tau_{s1}$)  the $\bZ_2$-gradings
on $S(TX, K_X) \otimes\otimes _{n=1}^\infty {\rm Sym} (TX_n) \otimes F^k_V
\otimes Q(W)$ defined by 
\begin{eqnarray}\label{f7}
\tau_{e1}= \tau_e \otimes 1 + 1 \otimes \tau_1,\qquad
\tau_{s1}= \tau_s \otimes 1 + 1 \otimes \tau_1.
\end{eqnarray}

Let $h^{V_v}$ be the metric on $V_v$ induced by the  metric $h^V$ on $V$.
In the following, we  identify $\Lambda V_v$ with 
$\Lambda \overline{V}^*_v$ by using the Hermitian metric $h^{V_v}$ on $V_v$.
By (\ref{hyp3}), 
as in (\ref{dirac4}), there is a natural isomorphism between 
$\bZ_2$-graded $C(V)$-Clifford modules over $F$,
\begin{eqnarray}\label{f6}
S(V_0^{\bf R}, {\otimes}_{0<v} (\det V_v)^{-1}) \widehat{\otimes}_{0<v} \Lambda 
V_v \simeq S(V)_{|F}.
\end{eqnarray}

By using the above notations, we rewrite (\ref{hyp8}),  on the fixed point set $F$,
for $p \in \bN$,
\begin{eqnarray}\label{f1}\begin{array}{l}
{\cal F}_p (X) = \bigotimes_{0<v} \Big (\bigotimes_{n=1}^\infty 
{\rm Sym} (N_{v,n}) 
\bigotimes_{\stackrel{n\in \bN,}{ n>pv}} {\rm Sym} (\overline{N}_{v,n}) \Big ) 
\bigotimes_{n=1}^\infty {\rm Sym} (TY_n),\\
 {\cal F}'_p (X) = \bigotimes_{ \stackrel{0<v, n\in \bN,}{0\leq n \leq pv}}
\Big (  {\rm Sym} (N_{v,-n})  \otimes \det N_v \Big ),\\
{\cal F}^{-p}(X) = {\cal F}_p (X) \otimes {\cal F}'_p (X).
\end{array}\end{eqnarray}

 Let $V_0= V_0^\sR \otimes _\sR \bC$. From  (\ref{hyp2}),  (\ref{f6}), we get
\begin{eqnarray}\label{f2}\begin{array}{l}
{\cal F}^0(X) = \bigotimes_{n=1}^\infty {\rm Sym} \Big (\oplus_{0<v}( N_{v,n} 
\oplus \overline{N}_{v,n} ) \Big ) \bigotimes_{n=1}^\infty {\rm Sym} (TY_n) \\
\hspace*{20mm} \bigotimes  {\rm Sym} (\oplus_{0<v} N_{v,0}) \otimes 
\det(\oplus_{0<v} N_{v}),\\
F^1_V =   \bigotimes_{n=1}^\infty \Lambda (\oplus_{0<v}( V_{v,n} \oplus 
\overline{V}_{v,n} ) \oplus V_{0,n}) \\
\hspace*{10mm} \otimes S(V_0^{\sR}, \otimes_{0<v} (\det V_v)^{-1})
 \otimes _{0<v}\Lambda ( V_{v,0}),\\
F^2_V =   \bigotimes_{0<n\in \sZ + 1/2} \Lambda (\oplus_{0<v}( V_{v,n} \oplus 
\overline{V}_{v,n}) \oplus V_{0,n}) ,\\
Q(W)=  \bigotimes_{n=0}^\infty \Lambda (\oplus_v \overline{W}_{v,n})
\bigotimes_{n=1}^\infty \Lambda (\oplus_v W_{v,n}).
\end{array}\end{eqnarray}

Now  we can reformulate Theorem \ref{Th24} as follows.
\begin{thm}
 For each $\alpha$,  $h,\ p\in \bZ$, $p>0$, 
$m\in {1\over 2} \bZ$,  for  $i=1,\ 2$, $\tau = \tau_{e1}$ or $\tau_{s1}$, 
 we have the following identity in $K_{G_y}(B)$,
\begin{eqnarray}\label{f3}
\begin{array}{l}
{\rm Ind}_\tau (D^{Y_\alpha}\otimes (K_W \otimes K_X^{-1})^{1/2}
 \otimes {\cal F}^{-p} (X)  \otimes F^i_V \otimes Q(W),\\
\hspace*{30mm} m + {1\over 2} p^2 e(N) +{1\over 2} p d'(N), h)\\
\hspace*{5mm} =(-1)^{pd'(W)}{\rm Ind}_\tau (D^{Y_\alpha} 
\otimes (K_W \otimes K_X^{-1})^{1/2} 
\otimes {\cal F}^0 (X)  \otimes F^i_V\\
\hspace*{30mm} \otimes Q(W)\otimes L^{-p},
m + ph + p^2 e, h ).
\end{array}\end{eqnarray}
\end{thm}

{\em Proof} : The rest of this section  is devoted to  a proof of  Theorem 3.1. 
\hfill $\blacksquare$

\subsection{ \normalsize Proof of Theorem 3.1}

Inspired by [{\bf T}, \S 7], as in [{\bf LiuMaZ}, \S 3], for $p\in \bN^*$, we define the shift operators, 
\begin{eqnarray}\begin{array}{l}
r_*: N_{v, n} \to N_{v, n+pv}, \qquad 
r_*: \overline{N}_{v, n} \to \overline{N}_{v, n-pv}, \\
r_*: W_{v, n} \to W_{v, n+pv}, \qquad 
r_*: \overline{W}_{v, n} \to \overline{W}_{v, n-pv}, \\
r_*: V_{v, n} \to V_{v, n+pv}, \qquad 
r_*: \overline{V}_{v, n} \to \overline{V}_{v, n-pv}. 
\end{array}\end{eqnarray}

Recall that $L(N),\ L(W),\ L(V)$ are the complex line bundles over $F$ 
defined by (\ref{hyp5}).
Recall also  that $L= L(N)^{-1}\otimes L(W) \otimes L(V)$ 
is a trivial complex line bundle 
over $F$, and $g\in S^1$ acts on it by multiplication by $g^{2e}$.

\begin{prop} For $p\in \bZ$, $p>0$, $i=1,\ 2$, 
 there are  natural isomorphisms of vector bundles over $F$,
\begin{eqnarray}\label{shift1}
\begin{array}{l}
r_* ({\cal F}^{-p} (X)) \simeq {\cal F}^{0} (X) \otimes L(N)^p,\\
r_* (F^i_V) \simeq F^i_V \otimes L(V)^{-p}.
\end{array}\end{eqnarray}
For any  $p\in \bZ$, $p>0$, there is a natural $G_y \times S^1$-equivariant isomorpism
 of vector bundles over $F$,
\begin{eqnarray}\label{shift2}
r_* (Q(W)) \simeq Q(W) \otimes L(W)^{-p}.
\end{eqnarray}
\end{prop}

{\em Proof} : The  equation  (\ref{shift1}) was proved in [{\bf LiuMaZ}, Prop. 3.1].
To prove  (\ref{shift2}), we only need to consider the shift operator on 
the following elements,
\begin{eqnarray}\label{FV1}
\begin{array}{l}
 Q_W =   \bigotimes_{n=0}^\infty \Lambda 
(\oplus_{v\neq 0} \overline{W}_{v,n} )
\bigotimes_{n=1}^\infty \Lambda (\oplus_{v\neq 0} W_{v,n}  ) .
\end{array}\end{eqnarray}
We compute easily that
\begin{eqnarray}\label{FV2}
\begin{array}{l}
r_* Q_W  =   \bigotimes_{n=0}^\infty \Lambda 
(\oplus_{v\neq 0} \overline{W}_{v,n-pv} )
\bigotimes_{n=1}^\infty \Lambda (\oplus_{v\neq 0} W_{v,n+pv}  ) .
\end{array}\end{eqnarray}
Let $h^W$ be a Hermitian metric on $W$. 
Let $h^{W_v}$ be the metric on $W_v$ induced by $h^W$.
As in [{\bf LiuMaZ}, \S 3], the
hermitian metric $h^{W_v}$ on $W_v$ induces a natural isomorphism 
of complex  vector bundles over $F$,
\begin{eqnarray}\label{isomV}
\Lambda ^i \overline{W}_{v} \simeq \Lambda ^{\dim W_v -i} W_{v} 
\otimes \det  \overline{W}_v. 
\end{eqnarray}

$\bullet$ If $v>0$, for $n\in \bN,\ 0\leq   n < pv$, $0\leq i \leq \dim W_v$, 
(\ref{isomV}) induces a natural $G_y \times S^1$-equivariant 
 isomorphism of complex vector bundles
\begin{eqnarray}\label{FV3}
\begin{array}{l}
\Lambda ^i \overline{W}_{v,n-pv} \simeq \Lambda ^{\dim W_v -i} W_{v, -n+pv} 
\otimes \det  \overline{W}_v. 
\end{array}\end{eqnarray}

$\bullet$ If $v<0$, for $n\in \bN,\ 0<  n \leq -pv$, $0\leq i \leq \dim W_v$, 
(\ref{isomV}) induces a natural $G_y \times S^1$-equivariant 
 isomorphism of complex vector bundles
\begin{eqnarray}\label{FV4}
\begin{array}{l}
\Lambda ^i {W}_{v,n+pv} \simeq \Lambda ^{\dim W_v -i} \overline{W}_{v, -n-pv} 
\otimes (\det \overline{W} _v)^{-1}. 
\end{array}\end{eqnarray}
{}From  (\ref{hyp5}), (\ref{FV3}) and (\ref{FV4}), we have  
\begin{eqnarray}\label{FV5}
\begin{array}{l}
\displaystyle{\bigotimes_{\stackrel{n\in \bN, v>0,}{0 \leq n < p v}}
\Lambda ^{i_n} \overline{W}_{v,n-pv} 
\bigotimes_{\stackrel{n\in \bN, v<0,}{0 < n \leq -p v}}
\Lambda ^{i'_n} {W}_{v,n+pv} }\\
\displaystyle{\hspace*{10mm}
\simeq \bigotimes_{\stackrel{n\in \bN, v>0,}{0 \leq n < p v}}
\Lambda ^{\dim W_v -i_n} {W}_{v, -n+pv} 
\bigotimes_{\stackrel{n\in \bN, v<0,}{0 < n \leq -p v}}
\Lambda ^{\dim W_v -i'_n} \overline{W}_{v,-n-pv}  \otimes L(W)^{-p}. }
\end{array}\end{eqnarray}
{}From  (\ref{FV2}), (\ref{FV5}), we get (\ref{shift2}).

The proof of Proposition 3.1 is complete. \hfill $\blacksquare$

\begin{prop} For $ p\in \bZ$, $p>0$,  $i=1,\ 2$, 
the $G_y$-equivariant bundle isomorphism induced by (\ref{shift1}) and 
 (\ref{shift2}), 
\begin{eqnarray}\label{tran0}\begin{array}{l}
r_* : S(TY,K_X \otimes_{0<v} (\det N_v)^{-1}) 
\otimes (K_W\otimes K_X^{-1})^{1/2}  \\
\hspace*{25mm}\otimes {\cal F}^{-p} (X) 
\otimes F^i_V  \otimes Q(W)   \\
\hspace*{10mm}\to S(TY, K_X \otimes_{0<v} (\det N_v)^{-1})
\otimes (K_W\otimes K_X^{-1})^{1/2} \\
\hspace*{25mm}\otimes {\cal F}^{0} (X) \otimes
F^i_V \otimes Q(W)  \otimes L^{-p},
\end{array}\end{eqnarray}
verifies the following identities
\begin{eqnarray}\label{tran1}
\begin{array}{l}
r_*^{-1}\cdot \bJ_H\cdot r_* =   \bJ_H,\\
r_*^{-1}\cdot P\cdot r_* =P +  p\bJ_H 
+p^2 e - {1 \over 2} p^2  e(N) - {p \over 2} d'(N).
\end{array}\end{eqnarray}
For the $\bZ_2$-gradings, we have
\begin{eqnarray}\label{tran2}
\begin{array}{l}
r_*^{-1} \tau_e r_* =  \tau_e,\qquad
r_*^{-1}\tau_s r_* =  \tau_s,\\ 
r_*^{-1} \tau_1 r_*= (-1)^{p d'(W)} \tau_1.
\end{array}\end{eqnarray}
\end{prop}

{\em Proof} : We divide the argument into several steps.

1) The first 
equation of (\ref{tran1}) is obvious.

2) a) From  [{\bf LiuMaZ}, (3.23)] and (\ref{hyp4}), for $i=1,2$, 
on $F^i_V$, we have  
\begin{eqnarray}\label{tran3}
r_*^{-1} P r_* = P+ p\bJ_H+ {1 \over 2}p^2 e(V). 
\end{eqnarray}

b) Note that on 
$\otimes_{0<v, 0\leq n \leq pv} \det N_v$, 
$\bJ_H$ acts as $p e(N) + d'(N)$.  
On $S(TY, K_X\otimes \det (\oplus_{0<v} N_v)^{-1})
\otimes (K_W\otimes K_X^{-1})^{1/2}$, $\bJ_H$ acts as 
$- {1 \over 2} d'(N)+{1\over 2} d'(W)$. From  (\ref{hyp4}), (\ref{f1}), 
on $S(TY, K_X\otimes \det (\oplus_{0<v} N_v)^{-1}) 
\otimes (K_W\otimes K_X^{-1})^{1/2} \otimes {\cal F}^{-p} (X)$,
\begin{eqnarray}\label{tran4}
r_*^{-1} P r_* =  P + p\bJ_H - p^2 e(N) -{1 \over 2}p (d'(N) +d'(W)).
\end{eqnarray}

c) From  (\ref{hyp4}), (\ref{FV5}), on 
$\bigotimes_{\stackrel{n\in \bN, v>0,}{0 \leq n < p v}}
\Lambda ^{i_n} \overline{W}_{v,n} 
\bigotimes_{\stackrel{n\in \bN, v<0,}{0 < n \leq -p v}}
\Lambda ^{i'_n} {W}_{v,n} $,
one has
\begin{eqnarray}\label{tran5}
\\
\begin{array}{l}
\displaystyle{r_*^{-1} P r_* =  \sum_{\stackrel{n\in \bN, v>0,}{0 \leq n < p v}}
(\dim W_v -i_n) (-n + pv) +
\sum_{\stackrel{n\in \bN, v<0,}{0 < n \leq -p v}}(\dim W_v -i'_n) (-n - pv) }\\
\displaystyle{\hspace*{5mm}=  P + p\bJ_H  + \sum_{\stackrel{n\in \bN, v>0,}{0 \leq n < p v}} 
(\dim W_v) (-n + pv) + \sum_{\stackrel{n\in \bN, v<0,}{0 < n \leq -p v}}
(\dim W_v) (-n - pv) }\\
\hspace*{5mm}=  P + p\bJ_H  +{1\over 2} p^2 e(W) + {1\over 2} p d'(W).
\end{array}\nonumber\end{eqnarray}

{}From  (\ref{hyp7}), (\ref{tran3}), (\ref{tran4}) and (\ref{tran5}), 
 we get the second equality of (\ref{tran1}).

3)  The first two identities of (\ref{tran2}) were proved in 
[{\bf LiuMaZ}, Proposition 3.2]. 

For the $\bZ_2$-grading $\tau_1$, it changes only on 
$\bigotimes_{\stackrel{n\in \bN, v>0,}{0 \leq n < p v}}
\Lambda ^{i_n} \overline{W}_{v,n} 
\bigotimes_{\stackrel{n\in \bN, v<0,}{0 < n \leq -p v}}
\Lambda ^{i'_n} {W}_{v,n} $.
{}From  (\ref{hyp4}), (\ref{FV5}), we get the last equality of (\ref{tran2}).

The proof of Proposition 3.2 is complete.
\hfill $\blacksquare$\\

{\em Proof of Theorem 3.1} : From (\ref{hyp7}), (\ref{f7}) and  Propositions  3.2,
we easily obtain Theorem 3.1.
\hfill $\blacksquare$\\

\newpage

\section{ \normalsize Proof of Theorem \ref{Th25}}
\setcounter{equation}{0}

In this section, we  prove Theorem \ref{Th25}. 
 As in [{\bf LiuMaZ}, \S 4], we will construct a family twisted 
Dirac operator on $M(n_j)$, the fixed point set of the
induced $\bZ_{n_j}$ action on $M$. 
 By applying our $K$-theory version of the equivariant family
index theorem to this operator,
 we prove Theorem \ref{Th25}.

This section is organized as follows: 
In Section 4.1, we construct a family Dirac operator on $M(n_j)$. 
In Section 4.2, by introducing a shift operator, we will relate both sides
 of equation (\ref{e8}) to the index bundle
  of the family Dirac operator 
on $M(n_j)$. In Section 4.3, we prove Theorem \ref{Th25}.

In this section, we  make the same assumptions and use the same 
notations as in  Sections 2, 3.

\subsection{ \normalsize The Spin$^c$ Dirac operator on $M(n_j)$}

Let $\pi: M\to B$ be a  fibration of compact manifolds with fiber $X$ 
and $\dim_{\sR} X= 2l$. We assume that  $S^1$ acts fiberwise on $M$,
and $TX$ has an $S^1$-invariant Spin$^c$ structure.
Let $F=\{F_\alpha \}$ be the fixed point set of the $S^1$-action on $M$. 
Then $\pi: F\to B$ is a fibration with compact fiber $Y$.
For $n\in \bN,\ n>0$, let $\bZ_n \subset S^1$ denote
the cyclic subgroup of order $n$. 

 Let $V$ be a real even dimensional 
vector bundle  over $M$ with an $S^1$-invariant spin structure.
Let $W$ be an $S^1$-equivariant complex vector bundle  over $M$.

For $ n_j\in \bN$, $n_j>0$,  
let $M(n_j)$ be the fixed point set of the
induced $\bZ_{n_j}$-action on $M$. 
Then $\pi: M(n_j)\to B$ is a fibration with compact fiber $X(n_j)$.
Let $N(n_j)\to M(n_j)$ be the  normal bundle to $M(n_j)$ in $M$. 
As in [{\bf LiuMaZ}, \S 4.1], we see that
$N(n_j)$ and $V$ can be decomposed, as    
real vector bundles over $M(n_j)$, to 
\begin{eqnarray} \label{dec2}\begin{array}{l}
\displaystyle{N(n_j) \simeq  \bigoplus_{0<v <n_j/2} N(n_j)_v \oplus N(n_j)_{{n_j \over 2}}^{\sR} ,}\\
\displaystyle{V|_{M(n_j)} \simeq V(n_j)_{0}^{\sR} \bigoplus_{0<v <n_j/2} V(n_j)_v 
\oplus V(n_j)_{{n_j \over 2}}^{\sR} }
\end{array}\end{eqnarray}
respectively. In (\ref{dec2}), the last term is understood to be zero 
when $n_j$ is odd. We  also denote  by $V(n_j)_{0}$,  $V(n_j)_{{n_j \over 2}}$,
$N(n_j)_{{n_j \over 2}}$ the corresponding complexification of 
the real vector bundles $V(n_j)_{0}^{\sR}$,  $V(n_j)_{{n_j \over 2}}^{\sR}$  
 and $N(n_j)_{{n_j \over 2}}^{\sR}$ on $M(n_j)$. 
Then $N(n_j)_v$, $V(n_j)_v$'s are  complex vector bundles over 
$M(n_j)$ with $g\in \bZ_{n_j}$
acting by $g^v$ on it. 

Similarly, we also have the following 
$\bZ_{n_j}$-equivariant decomposition of  $W$ on $M(n_j)$,
\begin{eqnarray} \label{dec3}
W= \oplus_{0\leq v <n_j} W(n_j)_v,
\end{eqnarray}
Here $W(n_j)_v$ is a  complex vector bundle over 
$M(n_j)$ with $g\in \bZ_{n_j}$
acting by $g^v$ on it.

It is essential for us to know that
 the vector bundles $TX(n_j)$ and $V(n_j)_{0}^{\sR}$ are orientable. 
For this we have
the following lemma which generalizes [{\bf BT}, Lemmas 9.4, 10.1], 

\begin{lemma}  Let $R$ be a real, even dimensional 
orientable vector bundle over a manifold $M$.
 Let $G$ be a compact Lie group. We assume that $G$ acts on $M$, 
and lifts to $R$. We assume that $R$ has a $G$-invariant 
Spin$^c$ structure. For $g\in G$, let $M^g$ be the fixed point set of $g$ 
on $M$. Let $R_0$ be the 
subbundle of $R$ over $M^g$ on which $g$ acts trivially.
Then $R_0$ is even dimensional and orientable.
\end{lemma}

{\em Proof} : Let $h^R$ be the metric on $R$ which is induced from 
the  Spin$^c$ structure on $R$. As $g$ preserves the Spin$^c$ structure of $R$,
 $g$ is an isometry on $R$ and preserves the orientation of $R$. 
On $M^g$, we have the following decomposition of real vector bundles,
\begin{eqnarray}
R= R_0 \oplus R_1. \nonumber
\end{eqnarray}
Since the only possible real eigenvalue of $g$ on $R_1$ is $-1$, 
and $\det (g_{|R_1})=1$ on $M^g$, we know that 
$\dim_\sR  R_1= \dim_\sR  R-\dim_\sR  R_0$ must be even. 
So $\dim_\sR R_0$ is even.

Let $K_R$ be the $G$-equivariant complex line bundle over $M$ which is induced 
by the Spin$^c$ structure of $R$.  Now the action of $g$ on the fiber 
of the complex spinor bundle $S(R, K_R)$ at $x\in M^g$ gives an element 
$\widetilde{g} \in {\rm Spin}^c(R_x)  = {\rm Spin} (R_{x}) \times _{\bZ_2} 
S^1 \subset C(R_x) \otimes_\sR \bC$,  here $C(R_x) $ is 
the Clifford algebra of $R_x$. Let $\widetilde{g}= (\widetilde{g}_1, s)$, with 
$\widetilde{g}_1 \in {\rm Spin} (R_{x})$, $s\in S^1$. 
Let $\rho: {\rm Spin} (R_{x}) \to SO(R_x)$ be the standard representation 
of ${\rm Spin} (R_{x})$, then $\rho (\widetilde{g}_1)= g$. So 
$\widetilde{g} c(a) = c(ga) \widetilde{g}$ for $a\in R_x$. 
Here we denote by $c(\cdot)$ the Clifford action. 
This means that $\widetilde{g}$ commutes $c(a)$ for $a\in R_{0x}$, 
so $\widetilde{g} \in {\rm Spin}^c(R_{1x}) = {\rm Spin} 
(R_{1x}) \times _{\bZ_2} 
S^1 \subset C(R_{1x})\otimes_\sR \bC$ 
and $\widetilde{g}_1 \in {\rm Spin} (R_{1x})$.

Let $e_1, \cdots, e_{2k}$ be an orthonormal basis of $R_{1x}$, 
then $e_{i_1}\cdots e_{i_j}$ $(1\leq i_1 < \cdots < i_j\leq 2k)$ 
is an orthonormal basis of the complex vector space 
$C(R_{1x})\otimes_\sR \bC$. We define 
$\sigma: C(R_{1x})\otimes_\sR \bC \to \det (R_{1x})\otimes_\sR \bC$ by 
\begin{eqnarray}\begin{array}{l}
\sigma(e_{i_1}\cdots e_{i_j})
= e_1 \wedge \cdots \wedge e_{2k} \quad {\rm if}  \quad j= 2k=\dim_\sR R_1,\\
\hspace*{25mm}=0 \quad {\rm otherwise}.
\end{array}\nonumber\end{eqnarray}
By [{\bf BGV}, Lemma 3.22], 
\begin{eqnarray}
| \sigma (\widetilde{g})| = | \sigma (\widetilde{g}_1)| 
= {\det}^{1/2} ((1-g_{|R_1})/2).
\end{eqnarray}
So $\sigma (\widetilde{g})$ is a nonvanishing section of 
$\det (R_1)\otimes_\sR \bC$, 
$\det (R_1)\otimes_\sR \bC$ is a trivial complex line bundle on $F$. 
So $\det (R_1)$ is trivial, and  $R_1$ is orientable.
So $R_0$ is orientable. 

This completes the proof of Lemma 4.1.
\hfill $\blacksquare$\\

By  Lemma 4.1,  $TX(n_j)$ and $V(n_j)_{0}^{\sR}$ are even dimensional and
orientable over $M(n_j)$. 
Thus $N(n_j)$ is   orientable over $M(n_j)$. By  (\ref{dec2}),
$N(n_j)_{{n_j \over 2}}^{\sR}$ and $V(n_j)_{{n_j \over 2}}^{\sR}$ 
are also even dimensional and orientable over $M(n_j)$. 
In the following, we fix  the orientations of 
 $N(n_j)_{{n_j \over 2}}^{\sR}$ and  $V(n_j)_{{n_j \over 2}}^{\sR}$ over $M(n_j)$.
We also fix  the orienations of $TX(n_j)$ and $ V(n_j)_{0}^{\sR}$ 
which are induced by (\ref{dec2}) and the orientations on $TX, V$,  
$N(n_j)_{{n_j \over 2}}^{\sR}$ and $V(n_j)_{{n_j \over 2}}^{\sR}$.

Let 
\begin{eqnarray} \label{rn1}
r(n_j)= {1\over 2} ( 1 + (-1)^{n_j}).
\end{eqnarray}

\begin{lemma} Assume that (\ref{hyp1})  holds.  Let 
\begin{eqnarray} \label{line1}\begin{array}{l}
L(n_j)= \bigotimes _{0<v <n_j/2} \Big ( \det (N(n_j)_v) \otimes 
\det (\overline {V(n_j)_v})\\
\hspace*{20mm} \otimes \det (\overline{W(n_j)}_v) 
\otimes \det (W(n_j)_{n_j-v} )\Big ) ^{(r(n_j)+1)v}
\end{array}\end{eqnarray}
be the complex line bundle over $M(n_j)$. Then we have

i)  $L(n_j)$ has an $n_j^{\rm th}$ root over $M(n_j)$. 

ii) Let 
  \begin{eqnarray}\label{line2}
\begin{array}{l}
L_1= K_X \bigotimes _{0<v <n_j/2} \Big ( \det (N(n_j)_v) \otimes 
\det (\overline {V(n_j)_v})\Big ) \\
\hspace*{40mm} \otimes \det (W(n_j)_{n_j/2}) \otimes L(n_j)^{r(n_j)/n_j},\\
L_2= K_X \bigotimes _{0<v <n_j/2} \Big ( \det (N(n_j)_v) \Big ) 
\otimes \det (W(n_j)_{n_j/2}) \otimes L(n_j)^{r(n_j)/n_j}.
\end{array}\end{eqnarray}
Let $U_1= T X(n_j)\oplus V(n_j)_{0}^{\sR}$ and 
$U_2 = TX(n_j) \oplus V(n_j)_{{n_j \over 2}}^{\sR}$.
Then  $U_1$ (resp. $U_2$) has a Spin$^c$ 
structure defined by  $L_1$ (resp. $L_2$).
\end{lemma}

{\em Proof} : Both statements follow from the proof of [{\bf BT}, Lemmas 11.3 and 11.4]. \hfill $\blacksquare$\\

Lemma 4.2 allows us, as we are going to see, to apply
the constructions and results in Section 1.1 to the fibration
$M(n_j)\rightarrow B$, which is the main concern of this section.

For $p_j \in \bN$, $p_j< n_j$, $(p_j, n_j)=1$, $\beta_j = {p_j \over n_j}$, let us write 
\begin{eqnarray}\begin{array}{l}
{\cal F} (\beta_j) =  \otimes _{0<n\in \bZ} {\rm Sym} (TX(n_j)_n) 
\bigotimes_{0<v <n_j/2} {\rm Sym} 
\Big ( \bigoplus_{0< n \in \bZ + p_j v/n_j} N(n_j)_{v,n}\\
\hspace*{20mm}\bigoplus_{0< n \in \bZ - p_j v/n_j} \overline{N(n_j)}_{v,n} \Big ) 
\otimes _{0<n\in \bZ + {1 \over 2}} {\rm Sym} (N(n_j)_{{n_j \over 2},n} ),
\end{array}\label{idFV1}\\
\begin{array}{l}
F^1_V(\beta_j)  = \Lambda \Big (\oplus_{0<n\in \bZ} V(n_j)_{0,n} \bigoplus_{0< v <n_j/2} 
\Big ( \bigoplus_{0<n \in \bZ + p_j v/n_j} V(n_j)_{v,n} \\
\hspace*{20mm}\bigoplus _{0< n \in \bZ - p_j v/n_j} \overline{V(n_j)}_{v,n} \Big ) 
\oplus _{0< n\in \bZ + {1 \over 2} } V(n_j)_{{n_j \over 2}, n} \Big ),\\
F^2_V(\beta_j) = \Lambda \Big (\oplus_{0<n\in \bZ} V(n_j)_{{n_j \over 2},n} 
\bigoplus_{0< v <n_j/2} 
\Big ( \bigoplus_{0<n \in \bZ + p_j v/n_j+ {1 \over 2}} V(n_j)_{v,n}    \\
\hspace*{20mm}\bigoplus _{0< n \in \bZ - p_j v/n_j +{1 \over 2}} 
\overline{V(n_j)}_{v,n}\Big ) \oplus _{0< n\in \bZ + {1 \over 2} } 
V(n_j)_{0, n} \Big ),\\
Q_W(\beta_j)= \Lambda \Big (\bigoplus_{0 \leq v < n_j} 
\Big ( \bigoplus_{0<n \in \bZ + p_j v/n_j} W(n_j)_{v,n} 
\bigoplus _{0\leq n \in \bZ - p_j v/n_j} \overline{W(n_j)}_{v,n} \Big ) \Big ).
\end{array}  \nonumber 
\end{eqnarray}
 
We denote by $D^{X(n_j)}$ the $S^1$-equivariant Spin$^c$-Dirac operator 
on $S(U_1,L_1)$  or $S(U_2, L_2)$ along the fiber $X(n_j)$ 
defined as in Section 1.1. 
We denote by $D^{X(n_j)} \otimes {\cal F}(\beta_j)
\otimes F^i_V(\beta_j)\otimes Q_W(\beta_j)$ $(i=1,\ 2)$ the corresponding twisted Spin$^c$ Dirac 
operator on $S(U_i,L_i) \otimes {\cal F}(\beta_j)
\otimes F^i_V(\beta_j)\otimes Q_W(\beta_j)$ along the fiber $X(n_j)$. \\

{\bf Remark 4.1.}  In fact, to define an $S^1$ (resp. $G_y$)-action on 
$ L(n_j)^{r(n_j)/n_j}$, 
one must replace the $S^1$-action by its $n_j$-fold action
(resp. the $G_y$-action by $G_{y^{1 / n_j}}$-action). Here by abusing 
notation, we still say an $S^1$ (resp. $G_y$)-action
 without causing any confusion.\\

In the rest of this subsection, we will reinterpret all of the above objects
when we restrict ourselves to $F$, the fixed point set of the $S^1$ action.
We will use the notation of Sections 1.1 and 2.

Let $N_{F/M(n_j)}$ be the normal bundle to $F$ in $M(n_j)$. Then by 
(\ref{hyp2}),
\begin{eqnarray}\label{idFV2}\begin{array}{l}
N_{F/M(n_j)}= \bigoplus _{0<v : v\in n_j \bZ} N_v,\\
TX(n_j) \otimes_\sR  \bC = TY\otimes_\sR \bC \oplus_{0<v, v\in n_j \bZ}
 (N_v \oplus \overline{N}_v).
\end{array}\end{eqnarray}
By (\ref{hyp2}), (\ref{hyp3}) and (\ref{dec2}), the restriction to $F$ of 
$N(n_j)_v $, $V(n_j)_v$ $(1\leq v \leq n_j/2)$ is given by 
\begin{eqnarray}\label{idFV3}\begin{array}{l}
\displaystyle{N(n_j)_v = \bigoplus _{0<v': v'=v\ \sm (n_j) } N_{v'} 
\bigoplus _{0<v': v'= -v\ \sm  (n_j) } \overline{N}_{v'},  }\\
\displaystyle{V(n_j)_v = \bigoplus _{0<v': v'=v\ \sm (n_j) } V_{v'} 
\bigoplus _{0<v': v'= -v\ \sm  (n_j) } \overline{V}_{v'}.}
\end{array}\end{eqnarray}
And 
\begin{eqnarray}\label{idFV4}
V(n_j)_0 = V_{0}^{\sR}  \otimes_\sR \bC  \bigoplus_{0<v, v=0\ \sm (n_j)}
 (V_v \oplus \overline{V}_v).
\end{eqnarray}
By (\ref{idFV2})-(\ref{idFV4}), 
we have the following identifications of  real vector bundles over $F$,
\begin{eqnarray}\label{idFV5}\begin{array}{l}
N(n_j)_{{n_j \over 2}}^{\sR} =  \bigoplus_{0<v, v={n_j \over 2}\ \sm (n_j)}
 {N}_v ,\\
TX(n_j) = TY \bigoplus_{0<v, v=0\ \sm (n_j)}
 {N}_v ,\\
V(n_j)_{0}^{\sR} = V_{0}^{\sR} \bigoplus_{0<v, v=0\ \sm (n_j)}
 V_v ,\\
V(n_j)_{{n_j \over 2}}^{\sR} = \bigoplus_{0<v, v={n_j \over 2}\ \sm (n_j)}
  V_v.
\end{array}\end{eqnarray}

By (\ref{hyp3}) and (\ref{dec3}), the restriction to $F$ of 
$W(n_j)_v$ $(0\leq v < n_j)$ is given by
\begin{eqnarray}\label{idFV6}
W(n_j)_v = \oplus_{v'=v \ \sm \ (n_j)} W_{v'}.
\end{eqnarray}

We denote by $V_0= V_{0}^{\sR} \otimes_\sR \bC$ the complexification of
 $V_{0}^{\sR}$ over $F$.  As $(p_j, n_j)=1$, we know that for $v\in\bZ$, 
$p_jv /n_j\in \bZ$ iff $v/n_j\in \bZ$. 
Also, $p_jv /n_j\in \bZ+ {1 \over 2}$ iff $v/n_j\in  \bZ+ {1 \over 2}$.
 From  (\ref{idFV2})-({\ref{idFV6}), we then get
\begin{eqnarray}\label{idFV7}\\
\begin{array}{l}
{\cal F} (\beta_j) = \otimes_{0<n\in \sZ} {\rm Sym} (TY_n)
\bigotimes_{0<v, v=0, {n_j \over 2}\ \sm (n_j) } 
\bigotimes_{0<n\in \bZ + {p_j v \over n_j}} 
{\rm Sym} (N_{v,n} \oplus \overline{N}_{v,n}) \\
\hspace*{10mm}  \bigotimes_{0<v'< n_j/2}
{\rm Sym} \Big ( \oplus_{v=v'\ \sm (n_j)} \Big ( \oplus_{0<n\in \bZ 
+ {p_j v \over n_j}} N_{v,n} 
\oplus_{0<n\in \bZ - {p_j v \over n_j}} \overline{N}_{v,n} \Big ) \\
\hspace*{15mm} \oplus_{v=-v'\ \sm (n_j)} \Big (  \oplus_{0<n\in \bZ 
+ {p_j v \over n_j}} N_{v,n} 
\oplus_{0<n\in \bZ - {p_j v \over n_j}} \overline{N}_{v,n} \Big ) \Big ),
\end{array}\nonumber\end{eqnarray}
\begin{eqnarray}\begin{array}{l}
F^1_V(\beta_j) = \Lambda  \Big [ \oplus_{0< n\in \bZ} V_{0,n} 
\bigoplus_{0<v, v=0, {n_j \over 2}\ \sm (n_j)} 
\Big ( \oplus_{0< n\in \bZ +  {p_j v \over n_j}} V_{v,n} 
 \oplus_{0< n\in \bZ -  {p_j v \over n_j}} \overline{V}_{v,n} \Big )\\
\hspace*{10mm}  \bigoplus_{0<v' < n_j/2} 
\Big ( \bigoplus_{v=v', -v'\ \sm (n_j)} 
\Big (  \oplus_{0< n\in \bZ +  {p_j v \over n_j}} V_{v,n} 
 \oplus_{0<n\in \bZ - {p_j v \over n_j}} \overline{V}_{v,n} 
\Big ) \Big )  \Big ],\\
F^2_V(\beta_j) = \Lambda \Big [\oplus_{0< n\in \bZ+ {1 \over 2}} V_{0,n} 
\oplus_{0<v, v=0, {n_j \over 2}\ \sm (n_j) } 
\Big ( \oplus_{0< n\in \bZ +  {p_j v \over n_j}+ {1 \over 2}} V_{v,n}  
 \oplus_{0< n\in \bZ -  {p_j v \over n_j}+ {1 \over 2}} \overline{V}_{v,n} 
\Big )\\
\hspace*{10mm}  \bigoplus_{0<v' < n_j/2} \Big ( \oplus_{v=v', -v'\ \sm (n_j)} 
\Big (  \oplus_{0< n\in \bZ +  {p_j v \over n_j}+ {1 \over 2}} V_{v,n}  
\oplus_{0<n\in \bZ - {p_j v \over n_j}+ {1 \over 2}} \overline{V}_{v,n} 
\Big ) \Big ) \Big ],\\
Q_W(\beta_j)= \Lambda \Big (\bigoplus_v 
\Big ( \bigoplus_{0<n \in \bZ + p_j v/n_j} W_{v,n} 
\bigoplus _{0\leq n \in \bZ - p_j v/n_j} \overline{W}_{v,n} \Big ) \Big ).
\end{array} \nonumber
\end{eqnarray}

Now, we want to compare the spinor bundles  over $F$.
{}From  (\ref{line1}), (\ref{line2}), (\ref{idFV3}) and  (\ref{idFV6}),  
we get that over $F$ we have the identities
\begin{eqnarray} \label{line3} \qquad  \begin{array}{l}
L(n_j)^{r(n_j) \over n_j} = \bigotimes_{0< v' < n_j/2} \Big ( 
\bigotimes _{v=v' \ \sm (n_j)} (\det N_v \otimes \det \overline{V}_v \otimes \det \overline{W}_v ) ^{2v'}   \\
\hspace*{25mm} \bigotimes _{v=-v'\ \sm (n_j)} 
(\det N_v \otimes \det \overline{V}_v\otimes \det \overline{W}_v) ^{-2v'} \Big )^{r(n_j)/n_j},\\
L_1=K_X \otimes L(n_j)^{r(n_j)/n_j} \bigotimes_{0< v' < n_j/2} \Big ( 
\bigotimes _{v=v' \ \sm (n_j)} (\det N_v \otimes \det \overline{V}_v)    \\
\hspace*{20mm}\bigotimes _{v=-v'\ \sm (n_j)} 
(\det N_v \otimes \det \overline{V}_v) ^{-1} \Big )
 \bigotimes_{v= {n_j \over 2} \ \sm  (n_j)} \det W_v ,\\
L_2= K_X \otimes L(n_j)^{r(n_j)/n_j}\bigotimes_{0< v' < n_j/2} \Big ( 
\bigotimes _{v=v'\  \sm (n_j)} \det N_v   \\
\hspace*{25mm} \bigotimes _{v=-v'\ \sm (n_j)} (\det N_v ) ^{-1} \Big )  
\bigotimes_{v= {n_j \over 2} \ \sm  (n_j)} \det W_v .
\end{array}\end{eqnarray}

{}From  (\ref{idFV5}), we have, over $F$,
\begin{eqnarray} \label{idFV8} \begin{array}{l}
TX(n_j) \oplus V(n_j)_{0}^{\sR} = TY \oplus V_{0}^{\sR} 
  \oplus _{0<v, v=0\ \sm (n_j)} (N_v \oplus V_v),\\
TX(n_j) \oplus V(n_j)_{{n_j \over 2}}^{\sR}  
= TY  \oplus _{0<v, v=0 \ \sm (n_j)} N_v 
\oplus_{0<v, v={n_j \over 2}\ \sm (n_j)}  V_v.
\end{array}\end{eqnarray}
Recall that the Spin$^c$ vector bundles $U_1$, $U_2$ have been defined in
Lemma 4.2. Denote by 
\begin{eqnarray} \label{idFV11} \\
\begin{array}{l}
\displaystyle{S(U_1,L_1)' = S \Big (TY\oplus V_{0}^{\sR}, 
L_1 \bigotimes_{\stackrel{0<v,}{v=0\ \sm (n_j)}} 
(\det N_v \otimes \det V_v)^{-1}  \Big ) 
\bigotimes _{\stackrel{0<v,}{v=0\ \sm (n_j)}} \Lambda V_v,  }\\
\displaystyle{S(U_2,L_2)' = S \Big (TY, L_2 \bigotimes_{\stackrel{0<v,}{v=0\ \sm (n_j)}} 
(\det N_v)^{-1} \bigotimes _{\stackrel{0<v,}{v={n_j \over 2}\ \sm (n_j)}}
( \det V_v)^{-1} \Big  ) \bigotimes _{\stackrel{0<v,}{v={n_j \over 2} 
\ \sm (n_j)}} \Lambda V_v.  }
\end{array}\nonumber
\end{eqnarray}
Then from (\ref{dirac4}) and 
(\ref{idFV11}), for $i=1,\ 2$, we have the following isomorphism of 
Clifford modules over $F$,
\begin{eqnarray} \label{idFV12} \begin{array}{l}
S(U_i,L_i) \simeq  S(U_i,L_i)' 
\otimes \Lambda (\oplus _{0<v, v=0\ \sm (n_j)} N_v).
\end{array}\end{eqnarray}
We define the $\bZ_2$ gradings on $S(U_i,L_i)'\ (i=1,\ 2)$ induced by the 
$\bZ_2$-gradings on $S(U_i,L_i)$ $(i=1,\ 2)$ and on 
$\Lambda (\oplus _{0<v, v=0 \ \sm (n_j)} N_v)$ such that the isomorphism
(\ref{idFV12}) preserves the $\bZ_2$-grading.

We introduce   formally the following complex line bundles over $F$,
\begin{eqnarray}   \qquad\begin{array}{l}
L'_1= \Big [ L_1^{-1} \otimes _{\stackrel{0<v,}{ v=0\ \sm (n_j)}} 
(\det N_v \otimes \det V_v)\otimes _{0<v} 
(\det N_v \otimes \det V_v)^{-1} \otimes K_X  \Big ] ^{1/2}, \\
L'_2= \Big [ L_2^{-1}  \otimes _{\stackrel{0<v,}{ v=0\ \sm (n_j)}} 
\det N_v \otimes _{\stackrel{0<v,}{ v=n_j/2\ \sm (n_j)}}  \det V_v
\otimes _{0< v} (\det N_v)^{-1}  \otimes K_X  \Big ]^{1/2}.
\end{array} \nonumber\end{eqnarray}
{}From   (\ref{dirac4}), Lemma 4.2 and the assumption that $V$ is
spin, one verifies easily that 
$c_1({L'_i}^2) = 0 \ \sm (2)$ for $i=1,\  2$.
Thus $L'_1,\ L'_2$ are well defined complex line bundles over $F$.
For the later use, we also write down the following expressions of
$L'_i$ ($i=1,\ 2$) which can be deduced from (\ref{line3}):
\begin{eqnarray}  \label{line4} \begin{array}{l}
L'_1= \Big [ L(n_j)^{-r (n_j)/n_j} \otimes _{ v={n_j \over 2}\ \sm (n_j)} 
(\det N_v \otimes \det \overline{V}_v \otimes \det \overline{W}_v) \Big ] ^{{1 \over 2}} \\
\hspace*{20mm}\otimes _{ 0< v \leq {n_j \over 2}\ \sm (n_j)} 
(\det N_v)^{-1} 
\otimes _{ {n_j \over 2}< v <n_j \ \sm (n_j)}  (\det V_v )^{-1},\\
L'_2= \Big [ L(n_j)^{-r (n_j)/n_j} \otimes _{ v={n_j \over 2}\ \sm (n_j)} 
( \det N_v  \otimes \det V_v \otimes \det \overline{W}_v) 
\Big ] ^{{1 \over 2}}\\
\hspace*{20mm}
\otimes _{ 0 < v \leq {n_j \over 2} \sm (n_j)} (\det N_v)^{-1}.
\end{array}\end{eqnarray}

{}From  (\ref{line3}), (\ref{idFV11}), and  the definition of 
$L'_i$ $(i=1,\ 2)$,
we get the following identifications of Clifford modules over $F$,
\begin{eqnarray} \label{idFV9}\qquad  \begin{array}{l}
S(U_1, L_1 )' \otimes L'_1 = S(TY, K_X \otimes _{0<v} (\det N_v)^{-1})
\otimes S(V_{0}^{\sR} , \otimes _{0<v} (\det V_v)^{-1})\\
\hspace*{30mm}\otimes \Lambda (\oplus _{0<v, v=0 \ \sm (n_j)} V_v),\\
S(U_2, L_2 )' \otimes L'_2 = S(TY, K_X\otimes _{0<v} (\det N_v)^{-1})
\otimes \Lambda (\oplus_{0<v, v={n_j \over 2} \ \sm (n_j)} V_v).
\end{array}\end{eqnarray}

Let
\begin{eqnarray}\label{vol5}\begin{array}{l}
\displaystyle{
\Delta (n_j,N) = \sum_{{n_j \over 2} < v' < n_j} 
\sum _{0<v=v'\ \sm (n_j)} \dim N_v 
+ o(N(n_j)_{{n_j \over 2}}^{\sR} ),  } \\
\displaystyle{\Delta (n_j,V) = \sum_{{n_j \over 2} < v' < n_j} 
\sum _{0<v=v'\ \sm (n_j)} \dim V_v 
+ o(V(n_j)_{{n_j \over 2}}^{\sR} ), }
\end{array}\end{eqnarray}
with $o(N(n_j)_{{n_j \over 2}}^{\sR} )=  0\ {\rm or}\  1$
 (resp. $o(V(n_j)_{{n_j \over 2}}^{\sR} )= 0\ {\rm or}\  1$), 
depending on whether the 
given orientation on $N(n_j)_{{n_j \over 2}}^{\sR} $ 
( resp. $V(n_j)_{{n_j \over 2}}^{\sR} $) agrees or disagrees 
with the complex orientation of $\oplus_{v={n_j \over 2}\ \sm (n_j)} N_v$
(resp. $\oplus_{v={n_j \over 2}\ \sm (n_j)} V_v$).

By [{\bf LiuMaZ}, \S 4.1], (\ref{idFV6}) and  (\ref{idFV12}),
 for the $\bZ_2$-gradings induced by $\tau_s$, the difference of the
$\bZ_2$-gradings of (\ref{idFV9}) is  $(-1)^{\Delta (n_j,N)}$; 
for the $\bZ_2$-gradings induced by $\tau_e$, the difference of the 
$\bZ_2$-gradings
of the first (resp. second) equation of (\ref{idFV9}) is 
$(-1)^{\Delta (n_j,N)+ \Delta (n_j,V)}$ 
(resp. $(-1)^{\Delta (n_j,N)
+o(V(n_j)_{{n_j \over 2}}^{\sR} ) }$).

\subsection{\normalsize The Shift operators}

Let $p\in {\bf N}^*$ be fixed. For any $1\leq j\leq J_0$,
inspired by [{\bf T}, \S 9], as in [{\bf LiuMaZ}, \S 4], 
we define the following shift operators $r_{j*}$:
\begin{eqnarray}\label{shif1}\begin{array}{l}
r_{j*}: N_{v,n} \to N_{v, n + (p-1) v + p_j v /n_j}, \quad 
r_{j*}: \overline{N}_{v,n} \to 
\overline{N}_{v, n - (p-1) v - p_j v /n_j}, \\
r_{j*}: W_{v,n} \to  W_{v, n + (p-1) v + p_j v /n_j}, \quad 
r_{j*}: \overline{ W}_{v,n} \to 
\overline{ W}_{v, n - (p-1) v - p_j v /n_j}, \\
r_{j*}: V_{v,n} \to V_{v, n + (p-1) v + p_j v /n_j}, \quad 
r_{j*}: \overline{V}_{v,n} \to 
\overline{V}_{v, n - (p-1) v - p_j v /n_j}.
\end{array}\end{eqnarray}

If $E$  is a combination of the above bundles, we denote by
$r_{j*} E$ the bundle on which the action of $P$ is changed in 
the above way. 

Recall that the vector bundles $F_V^i$ $(i=1,\ 2)$ have been defined in
(\ref{f2}).
{}From  (\ref{e4}), we see that 
\begin{eqnarray}\label{shif2}\begin{array}{l}
{\cal F}_{p,j}(X)= {\cal F}_{p}(X) \otimes {\cal F}'_{p-1}(X) 
\bigotimes_{(v,n)\in \cup_{i=1}^{j} I^p_i} \Big({\rm Sym} ( N_{v,-n}) 
\otimes  \det N_v \Big ) \\
\hspace*{30mm}\bigotimes_{(v,n)\in \overline{I}^p_{j}} 
{\rm Sym}( \overline{N}_{v,n}) .
\end{array}\end{eqnarray}

\begin{prop} There are  natural isomorphisms of vector bundles over $F$,
\begin{eqnarray}\label{shif3}
 \qquad\begin{array}{l}
r_{j*} {\cal F}_{p,j-1}(X) \simeq  {\cal F} (\beta_j) 
\bigotimes _{0<v, v=0\ \sm (n_j)} {\rm Sym} (\overline{N}_{v,0})\\
\hspace*{30mm} \otimes_{0<v} (\det N_v)^{[{p_j v \over n_j}] +(p-1)v +1} 
 \bigotimes_{0<v, v=0\ \sm (n_j)} (\det N_v)^{-1} ,\\

r_{j*}{\cal F}_{p,j}(X) \simeq {\cal F} (\beta_j) 
\bigotimes _{0<v, v=0\ \sm (n_j)} {\rm Sym} ({N}_{v,0})
\otimes_{0<v} (\det N_v)^{[{p_j v \over n_j}]+(p-1)v  +1},\\

r_{j*}F^1_V \simeq S(V_{0}^{\sR}, \otimes_{0<v} (\det V_v)^{-1}) 
\otimes F^1_V(\beta_j)    
\bigotimes _{0<v, v=0\ \sm (n_j)} \Lambda  (V_{v,0}) \\
\hspace*{30mm}
\otimes_{0<v} (\det \overline{V}_v)^{[{p_j v \over n_j}]+(p-1)v  } ,\\
r_{j*}F^2_V \simeq F^2_V(\beta_j) 
\bigotimes _{0<v, v={n_j \over 2}\ \sm (n_j)} \Lambda  (V_{v,0}) 
\otimes_{0<v} (\det \overline{V}_v)
^{[{p_j v \over n_j} + {1 \over 2}]+(p-1)v },\\
r_{j*} Q(W) \simeq Q_W(\beta_j) \otimes _{0<v} 
(\det \overline{W}_v )^{[{p_j v \over n_j}] +(p-1)v +1} 
 \bigotimes_{0<v, v=0\ \sm (n_j)} (\det \overline{W}_v)^{-1}\\
\hspace*{30mm}\otimes _{v<0} (\det {W}_v )^{[-{p_j v \over n_j}] -(p-1)v} .
\end{array} 
\end{eqnarray}
\end{prop}

{\em Proof} : The proof is similar to the proof of Proposition 3.1.

Note that, by (\ref{e2}), for $v \in J= \{ v\in \bN|$ There exists $\alpha$ 
such that $N_v\neq 0$ on $F_\alpha \}$, there are no integer in 
$]{p_{j-1} v \over n_{j-1}}, {p_j v \over n_j}[$. 
So for $v\in J$, the elements $(v,n)\in \cup_{i=1}^{i_0} I_i^p$ 
are $(v, (p-1)v +1)$, $\cdots, (v, (p-1)v + [{p_{i_0} v \over n_{i_0}}])$ 
for $i_0= j-1, \ j$. Furthermore,
\begin{eqnarray}\label{shif5}\begin{array}{l}
\displaystyle{[{p_{j-1} v \over n_{j-1}}]= [ {p_j v \over n_j}] -1
\quad  {\rm if} \quad   v = 0  \quad  \sm (n_j), }\\
\displaystyle{[{p_{j-1} v \over n_{j-1}}]= [ {p_j v \over n_j}] 
\quad  {\rm if} \quad   v \neq 0  \quad  \sm (n_j). }
\end{array}\end{eqnarray}
By using (\ref{f2}), (\ref{shif1}), (\ref{shif2}), (\ref{shif5}), we can prove the first four equalities of (\ref{shif3}) as in the proof of [{\bf LiuMaZ}, Proposition 4.1].

{}From  (\ref{isomV}), we have the natural  $G_y \times S^1$-equivariant 
isomorphisms of complex  vector bundles over $F$,
\begin{eqnarray}\label{shif6}\qquad \begin{array}{l}
\displaystyle{\bigotimes_{\stackrel{n\in \bN, v>0,}
{0 \leq n < (p-1) v + { p_j v \over n_j}}}
\Lambda ^{i_n} \overline{W}_{v,n-(p-1) v-{ p_j v \over n_j}} \simeq \bigotimes_{\stackrel{n\in \bN, v>0,}{0 \leq n < (p-1) v + { p_j v \over n_j}}}
\Lambda ^{\dim W_v -i_n} {W}_{v, -n+(p-1) v + { p_j v \over n_j}} }\\
\displaystyle{\hspace*{25mm}\bigotimes _{0<v} 
(\det \overline{W}_v )^{[{p_j v \over n_j}] +(p-1)v +1} 
 \bigotimes_{0<v, v=0\ \sm (n_j)} (\det \overline{W}_v)^{-1},}  \\
\displaystyle{
\bigotimes_{\stackrel{n\in \bN, v<0,}{0 < n \leq -(p-1) v-{ p_j v \over n_j}}}
\Lambda ^{i_n} {W}_{v,n+(p-1) v + { p_j v \over n_j}} }\simeq 
\bigotimes_{\stackrel{n\in \bN, v<0,}{0 < n \leq -(p-1) v-{ p_j v \over n_j}}}
\Lambda ^{\dim W_v -i_n} \overline{W}_{v,-n-(p-1) v-{ p_j v \over n_j}} \\
\displaystyle{\hspace*{50mm}
\bigotimes _{v<0} (\det {W}_v )^{[-{p_j v \over n_j}] -(p-1)v}. }
\end{array}\end{eqnarray}
{}From (\ref{f2}), (\ref{idFV7}), (\ref{shif6}), we get the last equation of (\ref{shif3}).

The proof of Proposition 4.1 is complete.
\hfill $\blacksquare$

\begin{lemma} Let us write 
\begin{eqnarray}\label{line5}\begin{array}{l}
L(\beta_j)_1 = L'_1 \otimes_{0<v} (\det N_v )^{[{p_j v \over n_j}]+(p-1)v  +1}
\otimes_{0<v}
(\det \overline{ V}_v)^{[{p_j v \over n_j}]+(p-1)v}    \\
\hspace*{20mm}
\otimes_{0<v, v=0 \ \sm (n_j)} (\det N_v )^{-1} \otimes _{v<0} (\det {W}_v )^{[-{p_j v \over n_j}] -(p-1)v}\\
\hspace*{20mm}\otimes _{0<v} 
(\det \overline{W}_v )^{[{p_j v \over n_j}] +(p-1)v +1} 
 \bigotimes_{0<v, v=0\ \sm (n_j)} (\det \overline{W}_v)^{-1},\\
L(\beta_j)_2 = L'_2 \otimes_{0<v} (\det N_v)^{[{p_j v \over n_j}] +(p-1)v +1}  
\otimes_{0<v} (\det \overline{ V}_v)
^{[{p_j v \over n_j} + {1 \over 2}]+(p-1)v  } \\
\hspace*{20mm}  \otimes_{0<v, v=0\ \sm (n_j)} (\det N_v )^{-1} \otimes _{v<0} (\det {W}_v )^{[-{p_j v \over n_j}] -(p-1)v} \\
\hspace*{20mm}\otimes _{0<v} 
(\det \overline{W}_v )^{[{p_j v \over n_j}] +(p-1)v +1} 
 \bigotimes_{0<v, v=0\ \sm (n_j)} (\det \overline{W}_v)^{-1}.
\end{array}\end{eqnarray}
Then $L(\beta_j)_1,\ L(\beta_j)_2$ can be  extended naturally to  
$G_y\times S^1$-equivariant
complex line bundles which we will still denote by 
$L(\beta_j)_1,\ L(\beta_j)_2$ respectively over $M(n_j)$.
\end{lemma}

{\em Proof} : Write 
\begin{eqnarray} \label{numb1}
[{p_j v \over n_j}] = {p_j v \over n_j} - {\omega(v)  \over n_j}.
\end{eqnarray}
Note that for $v={n_j \over 2}\ \sm (n_j)$, 
${\omega(v)  \over n_j} = {1 \over 2}$. 

We introduce the following line bundle over $M(n_j)$,
\begin{eqnarray}\label{line6}\begin{array}{l}
L^{\omega} (\beta_j)= \bigotimes_{0< v< {n_j \over 2}}
 \Big  (\det (N(n_j)_{v}) \otimes \det (\overline{V(n_j)}_{v}) \\
\hspace*{25mm} \otimes \det (\overline{W(n_j)}_v) 
\otimes \det (W(n_j)_{n_j-v} )\Big ) 
^{-\omega(v)- r(n_j)v} .
\end{array}\end{eqnarray}
 As in [{\bf LiuMaZ}, (4.38)], Lemma 4.2 implies 
$L^{\omega} (\beta_j)^{1/n_j}$ is well defined over $M(n_j)$.

The contributions of $N$ and $V$ in $L(\beta_j)_1, L(\beta_j)_2$ are the same 
as given in [{\bf LiuMaZ}, Lemma 4.2], we only need to calculate 
the contribution of $W$ in $L(\beta_j)_1, L(\beta_j)_2$.  
Actually from [{\bf LiuMaZ}, (4.37), (4.44)], 
(\ref{hyp5}),  (\ref{idFV6}), (\ref{line4}), (\ref{line5}), (\ref{numb1}) and 
(\ref{line6}), we get 
\begin{eqnarray}\begin{array}{l}
\displaystyle{
L(\beta_j)_1= L^{-(p-1)-p_j /n_j} \otimes L^{\omega} (\beta_j)^{1/n_j}
\bigotimes_{0< v \leq {n_j \over 2}} \det (\overline{W(n_j)}_v), }\\
\displaystyle{L(\beta_j)_2 = L^{-(p-1)-{p_j \over n_j}} 
\otimes L^{\omega}(\beta_j)^{{1 \over n_j}} 
\bigotimes_{0< v \leq  {n_j \over 2}} \det (\overline{W(n_j)}_v)}\\
\displaystyle{\hspace*{25mm}
\bigotimes_{1\leq m \leq p_j/2} \bigotimes_{m-{1 \over 2}<  p_j v' /n_j < m} 
\det (\overline{V(n_j)}_{v'}).}
\end{array}\end{eqnarray}

The proof of Lemma 4.3 is complete.\hfill $\blacksquare$\\

Let us write
\begin{eqnarray}\label{re1} \begin{array}{l}
\varepsilon (W)= -{1 \over 2} \sum_{0<v} (\dim W_v)\Big [ ( [{p_j v \over n_j}] + (p-1)v) ([{p_j v \over n_j}] + (p-1)v+1)\\
\hspace*{20mm}-( {p_j v \over n_j} + (p-1) v) 
\Big (2 \Big ([{p_j v \over n_j}] + (p-1) v \Big )+1 \Big ) \Big ]\\
\hspace*{10mm}-{1 \over 2} \sum_{v<0} (\dim W_v)\Big [ ( [-{p_j v \over n_j}]
 - (p-1)v) ([-{p_j v \over n_j}] - (p-1)v+1)\\
\hspace*{20mm}+( {p_j v \over n_j} + (p-1) v) 
\Big (2 \Big ([-{p_j v \over n_j}] - (p-1) v \Big )+1 \Big ) \Big ],
\end{array}\end{eqnarray}
\begin{eqnarray} \begin{array}{l}
\varepsilon_1 = {1 \over 2} \sum_{0<v} (\dim N_v - \dim V_v) \Big [ ( [{p_j v \over n_j}] + (p-1)v) ([{p_j v \over n_j}] + (p-1)v+1)\\
\hspace*{20mm}-( {p_j v \over n_j} + (p-1) v) 
\Big (2 \Big ([{p_j v \over n_j}] + (p-1) v \Big )+1 \Big ) \Big ],\\
\varepsilon_2 = {1 \over 2} \sum_{0<v} (\dim N_v)\Big [
([{p_j v \over n_j}] + (p-1)v) ([{p_j v \over n_j}] + (p-1)v+1)\\
\hspace*{20mm}-( {p_j v \over n_j} + (p-1) v) 
\Big (2([{p_j v \over n_j}] + (p-1) v)+1 \Big ) \Big ]\\
\hspace*{10mm}
-{1 \over 2}\sum_{0<v} (\dim V_v) \Big [ ( [{p_j v \over n_j}+ {1\over 2}]
 + (p-1)v)^2 \\
\hspace*{20mm}
-2 ( {p_j v \over n_j} + (p-1)v)([{p_j v \over n_j}+ {1\over 2}] 
+ (p-1)v) \Big ].
\end{array} \nonumber \end{eqnarray}
Then  $\varepsilon (W), \varepsilon_1,\ \varepsilon_2$ 
are locally constant functions on $F$.

Recall that the involutions $\tau_e, \tau_s$ and $\tau_1$ were defined
in Section 3.1. 
Also recall that if $E$ is a $S^1$-equivariant vector bundle over $M$, then 
the weight of the $S^1$-action 
on $\Gamma(F,E)$ is given by the action $\bJ_H$ (cf. \S 3.1).
\begin{prop}  For $i=1,\ 2$, the $G_y$-equivariant isomorphisms induced by 
(\ref{idFV9}) and (\ref{shif3}),
\begin{eqnarray} \label{rel0}\qquad \begin{array}{l}
r_{i1} : S(TY, K_X\otimes_{0<v} (\det N_v)^{-1})
\otimes (K_W\otimes K_X^{-1})^{1/2}  \\
\hspace*{20mm}\otimes {\cal F}_{p,j-1}(X) 
\otimes F^i_V \otimes Q(W) \to \\
\hspace*{10mm}S(U_i,L_i)' \otimes (K_W\otimes K_X^{-1})^{1/2}\otimes {\cal F}(\beta_j) \otimes F^i_V(\beta_j) \\
\hspace*{20mm}\otimes Q_W(\beta_j) \otimes L(\beta_j)_i
\otimes _{\stackrel{0<v,}{v=0 \sm (n_j)}} {\rm Sym} (\overline{N}_{v,0}),\\
r_{i2} : S(TY, K_X\otimes_{0<v} (\det N_v)^{-1}) 
\otimes (K_W\otimes K_X^{-1})^{1/2}\\
\hspace*{20mm} \otimes {\cal F}_{p,j}(X) 
\otimes F^i_V\otimes Q(W) \to \\
\hspace*{10mm}S(U_i,L_i)'\otimes (K_W\otimes K_X^{-1})^{1/2} \otimes {\cal F}(\beta_j) \otimes F^i_V(\beta_j) \\
\hspace*{20mm}\otimes Q_W(\beta_j) \otimes L(\beta_j)_i  
\otimes _{\stackrel{0<v,}{v=0\ \sm (n_j)}}( {\rm Sym} ({N}_{v,0}) 
\otimes \det N_v),
\end{array}\end{eqnarray}
have 
the following properties : 1) for $i=1,\ 2$, $\gamma=1,\ 2$, 
\begin{eqnarray}  \label{rel2} \begin{array}{l}
 \quad r_{i \gamma} ^{-1} \bJ_ H r_{i \gamma} = \bJ_H,\\
 \quad r_{i \gamma} ^{-1} P r_{i \gamma} = P + ({p_j \over n_j} 
+ (p-1)) \bJ_H + \varepsilon_{i \gamma},
\end{array}\end{eqnarray}
where 
\begin{eqnarray}\begin{array}{l}
\varepsilon_{i1} = \varepsilon_i+ \varepsilon (W) - e(p, \beta_{j-1}, N),\\
\varepsilon_{i2} = \varepsilon_i+ \varepsilon (W) - e(p, \beta_j, N).
\end{array}\end{eqnarray}

2) Recall that $o(V(n_j)_{n_j \over 2}^\sR)$ was defined in (\ref{vol5}). Let 
\begin{eqnarray}\label{rel3} \begin{array}{l}
\mu_1 = -\sum_{0<v} [{p_j v \over n_j}] \dim V_v + \Delta(n_j,N) 
+ \Delta(n_j, V) \quad \sm (2),\\
\mu_2 = - \sum_{0<v} [{p_j v \over n_j}+ {1 \over 2}] \dim V_v 
+ \Delta(n_j,N)+ o(V(n_j)_{n_j \over 2}^\sR) \quad \sm (2),\\
\mu_3 = \Delta(n_j,N)   \quad \sm (2)  ,\\
\mu_4= \sum_v (\dim W_v) ( [{p_j v \over n_j}] + (p-1)v) + \dim W
+ \dim W(n_j)_0 \quad \sm (2) .
\end{array}\end{eqnarray}
Then for $i=1,\ 2$; $\gamma= 1,\ 2$, 
\begin{eqnarray}\label{rel4}\begin{array}{l}
r_{i \gamma} ^{-1} \tau_e r_{i \gamma} = (-1)^{\mu_i} \tau_e, \quad
r_{i \gamma}^{-1} \tau_s r_{i \gamma}= (-1)^{\mu_3} \tau_s,\\
r_{i \gamma} ^{-1} \tau_1 r_{i \gamma}= (-1)^{\mu_4} \tau_1.
\end{array}\end{eqnarray}
\end{prop}

{\em Proof} : The first equality of (\ref{rel2}) is trivial. From (\ref{e6}) and (\ref{shif5}), one has
\begin{eqnarray}\label{rel10}
e(p,\beta_j, N)= e(p, \beta_{j-1}, N) 
+ \sum_{0<v, v=0\ \sm(n_j)}\Big ((p-1)v + {p_j v \over n_j} \Big ) \dim N_v.
\end{eqnarray}

Denote by $\varepsilon_i (V)\ (i=1,\ 2)$  the contribution of $\dim V$ in 
$ \varepsilon_i$ $(i=1,\ 2)$ respectively. 
Then from [{\bf LiuMaZ}, (4.52), (4.53)], on $F^i_V$, we have 
\begin{eqnarray}\label{rel11}\begin{array}{l}
 r_{j*}^{-1} P r_{j*}= P + ((p-1)+ {p_j \over n_j}) \bJ_H + \varepsilon_i (V) .
\end{array}\end{eqnarray}
{}From (\ref{shif6}), as in (\ref{tran5}), on $Q(W)$, we get
\begin{eqnarray}\label{rel12}
 r_{j*}^{-1} P r_{j*}= P + ((p-1)+ {p_j \over n_j}) \bJ_H 
+ \varepsilon (W) + {1 \over 2} \Big ( (p-1) + {p_j \over n_j} \Big ) d'(W).
\end{eqnarray}

{}From (\ref{rel10}), (\ref{rel11}), (\ref{rel12}), 
and by proceeding as in the proof of 
Proposition 3.2, as in [{\bf LiuMaZ}, Proposition 4.2], 
  one deduces easily the second equation of (\ref{rel2}).
 
Finally from  the discussion following (\ref{vol5}), 
and [{\bf LiuMaZ}, (4.50)],  we get the first two equations of (\ref{rel4}). 
 By (\ref{idFV6}) and (\ref{shif6}), we get the last equation of (\ref{rel4}).
 
The proof of Proposition 4.2 is complete.
\hfill $\blacksquare$\\

\begin{lemma} For each connected component $M'$ of $M(n_j)$, 
$\varepsilon_1+ \varepsilon (W)$, $\varepsilon_2+ \varepsilon (W)$ 
are independent on the connected component of $F$ in $M'$.
\end{lemma}

{\em Proof} : From  (\ref{hyp7}), (\ref{idFV4}), (\ref{idFV6}), 
 (\ref{numb1}) and (\ref{re1}),  we have   
\begin{eqnarray}\label{const1}\qquad \begin{array}{l}
\displaystyle{
\varepsilon_1= {1\over 2} \sum_{0\leq v' < n_j} \sum_{ v = v'\ \sm (n_j)}
(\dim N_v - \dim V_v-\dim W_v ) }\\
\displaystyle{\hspace*{30mm}\Big [ -( {p_j v \over n_j} + (p-1)v)^2 
- {\omega(v')(n_j -\omega( v')) \over n_j^2} \Big ] }\\
\displaystyle{\hspace*{5mm}= (p-1 + {p_j \over n_j})^2 e 
-{1 \over 16} \Big ( \dim_\sR N(n_j)_{n_j\over 2}^{ \sR} 
-\dim_\sR V(n_j)_{n_j\over 2}^{ \sR}
-2 \dim W(n_j)_{n_j\over 2}  \Big )   }\\
\displaystyle{\hspace*{10mm}- {1 \over 2} \sum_{0<v' < n_j/2} 
\Big ( \dim N(n_j)_{v'} -\dim V(n_j)_{v'} 
-\dim W(n_j)_{v'} }\\
\displaystyle{\hspace*{35mm}-\dim W(n_j)_{n_j-v'} \Big )
 { \omega(v')(n_j - \omega(v')) \over n_j^2}.}
\end{array}\end{eqnarray}
By (\ref{re1}), $\varepsilon_2-\varepsilon_1$ was given in [{\bf LiuMaZ}, (4.49)], it is independent on the connected 
component of $F$ in $M'$.

The proof of Lemma 4.4 is complete.\hfill $\blacksquare$\\

The following Lemma  was proved in [{\bf BT}, Lemma 9.3] and [{\bf T}, Lemma 9.6] 
(cf. [{\bf LiuMaZ}, Lemma 4.6]).

\begin{lemma} Let $M$ be a smooth  manifold on which $S^1$ acts. 
Let $M'$ be a connected   component  of $M(n_j)$, 
the fixed point set of the subgroup $\bZ_{n_j}$ of $S^1$ on $M$. 
Let $F$ be the fixed point set of the $S^1$-action on $M$. 
Let $V\to M$ be a real, oriented, even dimensional vector bundle to 
which the $S^1$-action on $M$ lifts. Assume that $V$ is Spin over $M$.
 Let $p_j \in ]0,n_j[,\  p_j\in \bN$ 
and $(p_j, n_j)= 1$, then
\begin{eqnarray}\label{const3}\begin{array}{l}
\sum_{0<v} (\dim V_v )[{p_j v \over n_j}] + \Delta (n_j, V) \quad \sm (2),\\
\sum_{0<v} (\dim V_v) [{p_j v \over n_j} 
+ {1 \over 2}]+ o(V(n_j)_{n_j/2}^{ \sR})\quad \sm (2) 
\end{array}\end{eqnarray}
are independent on the connected components of $F$ in $M'$.
\end{lemma}

Recall that the number $d'(p,\beta_j,N)$ has been defined in (\ref{e6}).

\begin{lemma} For each connected component $M'$ of $M(n_j)$,
$d'(p,\beta_j, N) + \mu_i + \mu_4 \ \sm (2) $ $(i=1,\ 2,\ 3)$ 
are independent on the connected component of $F$ in $M'$.
\end{lemma}

{\em Proof }: By (\ref{rel3}), and Lemma 4.5, 
 to prove Lemma 4.6, we only need to prove 
$$\sum_{0<v} (\dim N_v )([{p_j v \over n_j}]+(p-1)v) 
+ \Delta (n_j, N) + \mu_4 
\quad \sm (2)$$ 
 is independent on the connected components of $F$ in $M'$.
But by [{\bf BT}, Lemma 9.3], as $\omega_2(TX\oplus W)_{S^1}=0 $, we know that,
$\sm (2)$, 
\begin{eqnarray}\label{const5}
\sum_{0<v} (\dim N_v )[{p_j v \over n_j}] + \Delta (n_j, N) +
\sum_{v} (\dim W_v )[{p_j v \over n_j}]
 \end{eqnarray}
is independent on the connected components of $F$ in $M'$. From (\ref{e6}), (\ref{e13}),  (\ref{const5}), we get Lemma 4.6.

The proof of Lemma 4.6 is complete.
\hfill $\blacksquare$

\subsection{ \normalsize Proof of Theorem \ref{Th25}}

{}From (\ref{e6}), (\ref{idFV3}), (\ref{idFV6})  and (\ref{shif5}), 
we see that 
\begin{eqnarray}\label{const8} \qquad
\begin{array}{l}
\displaystyle{
\sum_{0<v} \dim N_v = \sum_{0<v < {n_j \over 2}} \dim N(n_j)_v + {1 \over 2 }
 \dim_\sR N(n_j)_{n_j/2}^{ \sR} + \sum_{0<v, v=0\ \sm (n_j)} \dim N_v,  }\\
\displaystyle{
d'(p,\beta_j, N) = d'(p,\beta_{j-1}, N) + \sum_{0<v, v=0\ \sm (n_j)} \dim N_v.}
\end{array} 
\end{eqnarray}

By Lemma 4.6,  (\ref{const8}), $d'(p,\beta_{j-1}, N) + \sum_{0<v} \dim N_v 
+ \mu_i + \mu_4 \  \sm (2)$ $(i=1,2,3)$
 are constant functions on each connected component of $M(n_j)$.

{}From  Lemma 4.3, one knows that the Dirac operator 
$D^{X(n_j)} \otimes F(\beta_j) \otimes F^i_V(\beta_j) \otimes Q_W(\beta_j)
\otimes L(\beta_j)_i $ $(i=1,\ 2)$ is well-defined on $M(n_j)$. Thus,
by using Proposition 4.2, Lemma 4.4,  (\ref{idFV12}) and (\ref{const8}), 
for $i=1,\ 2$, $h\in \bZ$, $m\in {1 \over 2} \bZ$,  $\tau= \tau_{e1}$ or $\tau_{s1}$,
and by applying both the first and the second
equations of Theorem 1.1 to each connected component of $M(n_j)$ separately, 
we get the following identity in $K_{G_y}(B)$,
\begin{eqnarray} \label{last}\\
\begin{array}{l}
\sum_\alpha (-1)^{d'(p,\beta_{j-1}, N) +\sum_{0<v} \dim N_v } 
{\rm Ind}_\tau (D^{Y_\alpha} \otimes (K_W\otimes K_X^{-1})^{1/2} 
\otimes {\cal F}_{p, j-1}(X) \\
\hspace*{30mm}\otimes  F^i_V\otimes Q(W), 
m +e(p, \beta_{j-1}, N), h)\\
= \sum_\beta (-1)^{d'(p,\beta_{j-1}, N) + \sum_{0<v} \dim N_v + \mu}
 {\rm Ind}_\tau (D^{X(n_j)} \otimes (K_W\otimes K_X^{-1})^{1/2} 
\otimes F(\beta_j) \\
\hspace*{20mm}\otimes F^i_V(\beta_j) \otimes Q_W(\beta_j)
\otimes L(\beta_j)_i, 
m+ \varepsilon_i+ \varepsilon (W) + ({p_j\over n_j} +(p-1)) h, h)\\
=\sum_\alpha (-1)^{d'(p,\beta_{j}, N)+\sum_{0<v} \dim N_v  } 
{\rm Ind}_\tau (D^{Y_\alpha} \otimes (K_W\otimes K_X^{-1})^{1/2} 
\otimes  {\cal F}_{p, j}(X)\\
\hspace*{30mm}\otimes  F^i_V\otimes Q(W), 
m+e(p, \beta_{j}, N), h).
\end{array} \nonumber
\end{eqnarray}
Here $\sum_\beta$ means the sum over all  connected components of $M(n_j)$.
 In (\ref{last}), if $\tau=\tau_{s1}$, then $\mu= \mu_3+\mu_4$; 
if $\tau= \tau_{e1}$, then $\mu = \mu_i+\mu_4$.

The proof of Theorem \ref{Th25} is complete.\hfill $\blacksquare$\\

\centerline{------------------------}
\vskip 6mm

Kefeng LIU,
Department of Mathematics, Stanford University, Stanford, CA 94305, USA.

{\em E-mail address}: kefeng@math.stanford.edu

\vskip 6mm

Xiaonan MA,
Humboldt-Universitat zu Berlin, Institut f\"ur Mathematik, unter den Linden 6,
D-10099 Berlin, Germany.

{\em E-mail address}: xiaonan@mathematik.hu-berlin.de

\vskip 6mm

Weiping ZHANG,
Nankai Institute of Mathematics, Nankai university,
Tianjin 300071, P. R. China.

{\em E-mail address}: weiping@nankai.edu.cn

\end{document}